\documentclass[12pt]{amsart}
\usepackage{amssymb}
\usepackage{amsmath}
\usepackage{versions}
\usepackage{color}
\numberwithin{equation}{section}
\newtheorem{theorem}{Theorem}[section]
\newtheorem{lemma}[theorem]{Lemma}
\newtheorem{corollary}[theorem]{Corollary}
\newtheorem{example}[theorem]{Example}

\newtheorem{definition}[theorem]{Definition}
\newtheorem{proposition}[theorem]{Proposition}

\newtheorem{remark}[theorem]{Remark}

\newcommand\beq{\begin{equation}}
\newcommand\eeq{\end{equation}}
\newcommand\re{\mathrm {Re~}}
\newcommand\im{\mathrm {Im~}}
\newcommand\ii{\mathrm i}
\newcommand\al{\alpha}

\newcommand\Ga{\Gamma}
\newcommand\de{\delta}
\newcommand\eps{\varepsilon}
\newcommand\la{\lambda}
\newcommand\up{\upsilon}
\newcommand\ups{\upsilon}
\newcommand{\ph}{\varphi}
\newcommand\si{\sigma}
\newcommand\ta{\theta}
\newcommand\D{\mathbb D}
\newcommand\T{\mathbb T}
\newcommand\C{\mathbb C}
\newcommand\G{\mathcal{G}}
\newcommand\R{\mathbb R}
\newcommand\e{\mathrm e}

\newcommand\E{\mathcal E}
\newcommand\B{\mathcal B}

\newcommand\half{{\tfrac{1}{2}}}
\newcommand\dia{\diamondsuit}

\newcommand\df{\stackrel{\rm def}{=}}
\newcommand\Aut{\mathrm{Aut}~}
\newcommand\Schur{\mathcal{S}}

\newcommand\bbm{\begin{bmatrix}}
\newcommand\ebm{\end{bmatrix}}
\newcommand\ds{\displaystyle}
\DeclareMathOperator\aut{Aut}

\DeclareMathOperator\hol{Hol}
\def\proof{\noindent {\bf Proof}. }
\def\qed{\hfill $\square$ \vspace{3mm}}
\let\phi\varphi
\numberwithin{equation}{section}

\begin{document}

\title{3-extremal holomorphic maps and the symmetrised bidisc}
\author{Jim Agler, Zinaida A. Lykova and N. J. Young}
\date{26th July 2013}

\begin{abstract} 

We analyse the $3$-extremal holomorphic maps from the unit disc $\D$ to the symmetrised bidisc $
\G \df \{(z+w,zw): z,w\in\D\}$ with a view to the complex geometry and function theory of  $\G$. These are the maps whose restriction to any triple of distinct points in $\D$ yields interpolation data that are only just solvable.   We find a large class of such maps; they are rational of degree at most $4$. It is shown that there are  two qualitatively different classes of rational $\G$-inner functions of degree at most $4$, to be called {\em aligned} and {\em caddywhompus} functions; the distinction relates to the cyclic ordering of certain associated points on the unit circle. The aligned ones are $3$-extremal. 
 We describe a method for the construction of aligned rational $\G$-inner functions; with the aid of this method we reduce the solution of a $3$-point interpolation problem for aligned holomorphic maps from $\D$ to $\G$ to a collection of classical Nevanlinna-Pick problems with mixed interior and boundary interpolation nodes. 
Proofs depend on a form of duality for $\G$.
\end{abstract}

\thanks{The first author was partially supported by National Science Foundation Grant on  Extending Hilbert Space Operators DMS 1068830. The third author was partially supported by the UK Engineering and Physical Sciences Research Council grants EP/J004545/1 and EP/K50340X/1.The collaboration was partially supported by London Mathematical Society Grant 41219.}

\maketitle
\markboth{Jim Agler, Zinaida A. Lykova and N. J. Young}{3-extremal holomorphic maps and the symmetrised bidisc}
\section*{Contents} \label{contents}

\ref{newintro}. Introduction \hfill  Page \pageref{newintro}

\ref{prelim}.  Preliminaries \hfill \pageref{prelim}

\ref{duality}.  A form of duality for the symmetrised bidisc  \hfill \pageref{duality}

\ref{extremal}.  Extremal solvability \hfill \pageref{extremal}


\ref{maintheorem}.  The main theorem \hfill \pageref{maintheorem}

\ref{Eclasses}.  The classes $\mathcal E_{\nu n}$ of rational functions \hfill \pageref{Eclasses}

\ref{calculation}. Calculation of interpolating functions 
\hfill \pageref{calculation}

\ref{property}. Properties of interpolating functions \hfill   \pageref{property}

\ref{cancell}. Cancellations in some rational functions  \hfill \pageref{cancell}

\ref{snlemma}.  Snares  \hfill \pageref{snlemma}

\ref{boundary}. A bound for $s$ \hfill\pageref{boundary}

\ref{mainproof}.   Proof of the main theorem \hfill\pageref{mainproof}

\ref{caddy}.  Caddywhompus functions  \hfill \pageref{caddy}

\ref{dataOnBoundary}.  Target data on the boundary \hfill\pageref{dataOnBoundary}

\ref{doesnot}.  Weak solvability does not imply solvability \hfill\pageref{doesnot}

\ref{extsolvdata}.  More about extremally solvable data\hfill\pageref{extsolvdata}

\ref{conclude}. Concluding reflections \hfill\pageref{conclude}

 References \hfill\pageref{bibliog}

\section{Introduction}\label{newintro}
Hyperbolic geometry in the sense of Kobayashi \cite{Ko98} studies a domain $\Omega$ by means of the embedding of holomorphic discs in $\Omega$.  That is, it makes use of the elements of $\hol(\D,\Omega)$, the space of holomorphic maps from the open unit disc $\D$ of the complex plane into $\Omega$.  Here we study the hyperbolic geometry of the {\em open symmetrised bidisc}
$$
\G \df \{(z+w, zw): |z|< 1, |w|< 1\} \;\; \text{in} \;\; \C^2,
$$
but whereas  the Kobayashi distance in a domain $\Omega$ is defined in terms of maps in $\hol(\D,\Omega)$ whose images pass through a given {\em pair} of points in $\Omega$,  this paper is concerned with holomorphic maps  from $\D$ to $\G$ passing through a given {\em triple} of points.  One could think of this more delicate issue as constituting a form of `Kobayashi curvature'; it also relates to questions of interpolation that arise in an intended application to $H^\infty$ control.

As in Kobayashi's theory, there will be an emphasis on extremality.  The $3$-extremal holomorphic maps of the title are maps in $\hol(\D,\G)$ whose restriction to any $3$-point set yields interpolation data that are only just solvable.  This notion was introduced in \cite{ALY12}.  Formally, for any domain $\Omega$, a map  $h\in \hol(\D,\Omega)$  is {\em $n$-extremal} if, for any choice of $n$ distinct points $\la_1,\dots,\la_n$ in $\D$ and  for any open neighbourhood $U$ of  the closed unit disc, there is no function $f\in\hol(U,\Omega)$ such that $f(\la_j)=h(\la_j)$ for $j=1,\dots,n$.

 $\G$ was first studied because of its connection with a problem in control engineering, but it has turned out that the geometry of $\G$ is also significant for the theory of invariant distances \cite {JP,costara,EZ,NiPfZw}.

The $2$-extremal maps in $\hol(\D,\G)$ are precisely the complex geodesics of $\G$.  They are rational functions of degree at most $2$ and can be written down explicitly \cite{AY06,PZ}.  These geodesics are also {\em a fortiori} $3$-extremal maps, but the class of $3$-extremals is much larger.  In this paper we identify a large class of $3$-extremal maps in $\hol(\D,\G)$; they are rational functions of degree at most $4$.  They are also {\em $\G$-inner}, which means that they map almost every point of the unit circle to the distinguished 
boundary $b\G$ of $\G$ (Definition \ref{defGainner}).   Now 
\[
b\G= \{(z+w,zw):  |z|=1=|w| \}  \subset \C^2
\]
which is topologically a M\"obius band.  The fact that the distinguished boundary of $\G$ (unlike that of the bidisc) itself has a boundary  lends an additional richness to the function theory of $\G$.  A consequence relevant to this paper is that there are two qualitatively different classes of rational $\G$-inner functions of degree at most $4$, which we call {\em aligned} and {\em caddywhompus}\footnote{From the Urban Dictionary: caddywhompus -  something that is all out of wack, crooked, off centered, or not lined up correctly}; the distinction relates to the cyclic ordering of points on the unit circle $\T$ at which the values of the function lie on the edge of the M\"obius band.   We prove that aligned rational $\G$-inner functions of degree at most $4$ are $3$-extremal.

The heart of the paper is a technique for constructing aligned rational $\G$-inner functions of degree at most $4$, and the crux of the proof is a technical lemma (the `Snare Lemma' in Section \ref{snlemma}) which enables us to prove an appropriate boundedness property.   The method depends on certain `magic functions' $\Phi_\omega$ on $\G$, where $|\omega|=1,$ which play a role analogous to linear functionals in linear duality theory.  This special form of duality for $\G$ is described in Section \ref{duality}.

Our main result gives necessary and sufficient conditions for a $3$-point interpolation problem to be solvable by an aligned rational $\G$-inner function of degree at most $4$, in the sense that it reduces the problem to a collection of one-variable interpolation problems each of which has a classical solvability criterion.
We state the theorem, though some of the terminology will only be explained later.  Problem $\dia$ (see page \pageref{def1p}) is a one-variable Nevanlinna-Pick-type interpolation problem, with both interior and boundary interpolation nodes. Condition $\mathcal{C}_1(\la, z)$ (Definition \ref{condC},  page \pageref{condC}) is a parametrised family of Pick conditions (that is, the positivity of a family of matrices).
\begin{theorem} \label{main}
Let $\lambda_1,\lambda_2, \lambda_3$ be distinct points in $\mathbb{D}$
and let $z_1,z_2,z_3 \in \G$. 
The following statements are equivalent.
\begin{enumerate}
\item  There exists an aligned $\G$-inner function $h$ of degree at most $ 4$ such that
$h(\lambda_j) = z_j$ for  $ j=1, 2, 3$;

\item  condition $\mathcal{C}_1(\la, z)$ holds extremally and actively,
  and the associated Problem $\dia$ is solvable.
\end{enumerate}
\end{theorem}
The proof of the theorem is constructive, so that when (2) holds, we can in principle construct the desired function $h$, which will necessarily be $3$-extremal.  Corollary \ref{mainbis} gives a criterion for condition (2) to hold in terms of the rank and positivity of an associated matrix.

The definition of $3$-extremality that we introduced in \cite{ALY12} is not the only natural one; 
several others are possible. A secondary theme of the paper is to find relations between these notions and to explore which of them are fruitful --  see especially Section \ref{extsolvdata}.

\section{Preliminaries}\label{prelim}

We shall denote by $\Delta$ the closed unit disc and by $\Schur$ the {\em Schur class}, that is, the set $\hol(\D,\Delta)$ of holomorphic maps from $\D$ to $\Delta$.  The Riemann sphere will be denoted  by $\C^*$. 

In addition to the symmetrised bidisc $\G$ we shall also need its closure $\Gamma$, that is, the {\em closed symmetrised bidisc}
\[
\Gamma\df\{(z+w,zw): |z|\leq 1, |w|\leq1\}.
\]
Points in $\Gamma$ or $\G$ will be denoted by the symbols $(s,p)$, chosen to suggest `sum' and `product'.   The degree of a rational function $f$ will be denoted by $d(f)$.

By the {\em finite interpolation problem} for a subset $E$ of $\C^N$ we shall mean\\

\noindent {\bf Problem $IE$}:  {\em Given $n$ distinct points $\la_1, \dots, \la_n$ in the open unit disc $\D$ and $n$ points $z_1, \dots , z_n$ in $E$, find if possible an analytic function $h:\D\to E$ such that $h(\la_j)=z_j$ for $j=1,\dots,n.$  In particular, find a criterion for the existence of such an $h$.}

Interpolation data
\[
\la_j\in\D \mapsto z_j\in E, \quad j=1,\dots,n,
\]
will be said to be {\em solvable} if there exists $h\in\hol(\D,E)$ such that $h(\la_j)=z_j$ for each $j$.

In order to understand $3$-extremal holomorphic maps one must be concerned with Problem $I\Gamma$ or $I\G$ with $n=3$,   but statements will be formulated for general $n$ and $E$ where possible.  In the case that some target point $z_j$ lies in the topological boundary $\partial\Gamma$ of $\Gamma$ then Problem $I\Gamma$ can be solved relatively easily: see Section \ref{dataOnBoundary}.  The paper is therefore mainly concerned with the case that the target points $z_j$ are all in $\G$. 

If a problem $I\Ga$ is solvable then it has a solution $h\in\hol(\D,\Ga)$ that is  $\G$-inner; we may therefore restrict ourselves to the search for $\G$-inner interpolating functions. If a certain conjecture that we made in \cite{ALY12} is true then {\em all} $3$-extremal maps in $  \hol(\D,\G) $ are rational $\G$-inner functions of degree at most $4$  \cite[Theorem 7.3 and Observation 9.2]{ALY12}.  The conjecture is given below as the statement \eqref{mainQ}.

\section{A form of duality for the symmetrised bidisc}\label{duality}
A fruitful theme in hyperbolic geometry is a duality between $\hol(\D,\Omega)$ and $\hol(\Omega,\D)$ that culminates in a theorem of Lempert to the effect that the Lempert function and Carath\'eodory distance coincide for convex domains $\Omega \subset \C^d$.  The meaning of the statement is that, for any pair of points $z_1,\ z_2\in\Omega$, the two quantities
\begin{align}\label{lempert}
\delta_\Omega (z_1,z_2) &= \inf \{\rho(\la_1,\la_2):   \mbox{ there exist } \la_1,\la_2\in\D \mbox{ and } \notag \\
	&\hspace*{1cm}  h\in\hol(\D,\Omega) \mbox{ such that } h(\la_j)=z_j, \ j=1,2\} ,
\end{align}
and
\begin{align}\label{carath}
C_\Omega(z_1,z_2)  & = \sup \{\rho(F(z_1),F(z_2)): \quad F\in\hol(\Omega,\D)\}
\end{align}
are equal,  where $\rho$ denotes the pseudohyperbolic distance on $\D$.  The functions $\delta_\Omega$ and $C_\Omega$ are defined to be the Lempert function and the Carath\'eodory distance on $\Omega$ respectively \cite{JP}.

Lempert's theorem, and the theory of invariant distances of which it is a high point, suggest that, for any pair of domains $D$ and $\Omega$,  we should associate with the interpolation data
\beq\label{genInterp}
\la_j\in D \mapsto z_j \in\Omega, \quad j=1,\dots,n,
\eeq
and any $g\in\hol(\Omega,\D)$ the derived interpolation problem
\beq\label{derived}
\la_j\in D \mapsto g(z_j) \in \D, \quad j=1,\dots,n.
\eeq
\begin{definition}\label{wkslvbl}
The interpolation data \eqref{genInterp} are said to be {\em weakly solvable} if, for every $g\in\hol(\Omega,\D)$, the interpolation data \eqref{derived} are solvable.
\end{definition}

Clearly solvable data are weakly solvable (if $h$ solves \eqref{genInterp} then $g\circ h$ solves \eqref{derived}).  In some cases the converse is also true: if $\Omega$ is a polydisc then weak solvability implies solvability, for one may let $g$ run through the co-ordinate functions.  However, $\hol (\Omega,\D)$ may be a very small set (for example, if $\Omega=\C$), and so one cannot expect weak solvability to imply solvability in the absence of suitable properties of $\Omega$.  Nevertheless, since $\Gamma$ is so closely related to the bidisc, one could hope that weak solvability might imply solvability for Problem $I\Gamma$.  Alas, it is not so: see Section \ref{doesnot} below.  Consequently weak solvability is inadequate for the purpose of solving Problem $I\Gamma$. We therefore introduce a stronger form of duality specific to $\G$.

To explain this duality we use some special functions in $\hol(\G,\D)$ which enjoy a certain extremality property.  These are the rational functions $\Phi(\omega,\cdot)$ for $\omega\in\T$, where, for $(z, s, p) \in \mathbb{C}^3$ such that $zs \neq 2$,
\[
\Phi(z,s,p)= \frac{2z p-s}{2-z s}.
\]
These functions have the property that 
\beq\label{Phicriter}
\left| \Phi(z,s,p)\right| < 1 \mbox{ for all } z\in\Delta \Leftrightarrow (s,p)\in\G
\eeq																									
(see, for example, \cite{AY04}).  In complex-geometric terms, the $\Phi(\omega,\cdot)$ are the {\em magic functions} \cite{magic} of $\G$  (though we shall not use this fact). 
The  function  $\Phi(z, s, p)$ 
is defined for $(z, s, p) \in \mathbb{C}^3$ such that $zs \neq 2$.  In
particular, $\Phi$ is defined and analytic on $\mathbb{D} \times
\Gamma$ (since $|s| \le 2$ when $(s, p) \in \Gamma$). We shall write $\Phi_z(s, p)$ as a synonym for $\Phi(z, s, p)$. The function $\Phi$ plays a central role in the study of $\Gamma$. See \cite{AY2} for an account of how $\Phi$ arises.

With the problem
\beq\label{primal}
\la_j\in\D \mapsto z_j\in\G, \quad j=1,\dots,n,
\eeq
and any function $m \in\Schur$, we associate the Nevanlinna-Pick problem
\beq\label{derivedNP}
\la_j\in\D \mapsto \Phi(m(\la_j), z_j) \in \D \quad j=1,\dots,n.
\eeq
  In this way a problem $I\Gamma$ is associated with a family of classical Nevanlinna-Pick problems.  It is the study of these associated problems that constitutes the stronger form of duality that we consider.

For $\alpha \in \C$ we write
$$
B_\alpha(z) = \frac{z-\alpha}{1-\overline \alpha z}.
$$
In the event that $\al\in\D$ the rational function $B_\al$ is called a {\em Blaschke factor}.
A {\em M\"obius function} is a function of the form $cB_\alpha$
for some $\alpha \in \mathbb{D}$ and $c\in \mathbb{T}$.
The set of all M\"obius functions is the automorphism group 
$\Aut \mathbb{D}$ of $\mathbb{D}$. 
We denote by $\B l_n$ the set of Blaschke products of degree at most $n$.
\begin{definition}\label{condC}
Interpolation data
\begin{equation}
\label{Gdata}
\la_j\in\D \mapsto z_j\in\G, \quad j=1,\dots,n,
\end{equation}
where $\lambda_1, \dots, \lambda_n$ are distinct points in $\D$,
{\em satisfy condition $\mathcal{C}$} (or {\em condition $\mathcal{C}_\nu$}, for some non-negative integer $\nu$) if, for every $\up\in \Schur$ (or for every $\up \in \B l_{\nu}$, respectively), the Nevanlinna-Pick data
\begin{equation}
\label{upsdata}
\la_j\in\D \mapsto \Phi(\up(\la_j),z_j)\in\D, \quad j=1,\dots,n,
\end{equation}
are solvable.
\end{definition}
Clearly, if $h\in\hol(\D,\G)$ is a solution of the problem \eqref{primal} then, for any $m\in\Schur$,  $\Phi\circ(m,h)$ is  a solution of the derived problem \eqref{derivedNP}.  Thus condition $\mathcal C$ and, {\em a fortiori}, condition $\mathcal C_\nu$ are necessary for the solvability of Problem $I\Gamma$.  The question arises:
is condition $\mathcal C$ sufficient for the solvability of Problem $I\Gamma$?  If this question (originally posed in \cite{ALY12}) can be answered affirmatively it will be a major step towards the numerical solution of the spectral Nevanlinna-Pick problem for $2\times 2$ matrix functions, a problem that is currently poorly understood.  In \cite{ALY12} we conjectured that $\mathcal C$ {\em is} sufficient for solvability, and indeed that the following stronger statement is true:
\begin{equation} \label{mainQ}
\textit {Condition $\mathcal{C}_{n-2}$ is sufficient for the solvability of Problem I$\Ga$. }
\eeq
   Theorem \ref{main} is a partial affirmative answer in the case $n=3$.

 An important role will be played by the analogue for  $\hol(\D,\Ga)$ of inner functions, defined as follows.
\begin{definition}\label{defGainner}
A function $h\in\hol(\D,\Ga)$ is {\em $\Ga$-inner} or {\em $\G$-inner} if
\[
\lim_{r\to 1-} h(r\la) \in b\Ga
\]
for almost all $\la\in\T$ with respect to Lebesgue measure, where $b\Ga$ denotes the distinguished boundary
$\{(z+w,zw): |z|=1=|w|\}$ of $\Ga$.
\end{definition}
As was mentioned in Section \ref{prelim}, if a problem $I\Ga$ is solvable then it has a $\G$-inner solution, so there is no loss in seeking $\G$-inner interpolating functions.

\section{Extremal solvability}\label{extremal}

In classical Nevanlinna-Pick theory interpolation data that are {\em extremally} solvable  admit a {unique} interpolating function $q$, which is a Blaschke product of degree less than the number of interpolation nodes.   Moreover, there is a simple formula for $q$ in terms of a null vector of the Pick matrix of the data (for example, \cite{AgMcC}).  Extremally solvable data play an important role in the present study too.  In this section we introduce a natural geometric notion of extremally solvable interpolation data, as well as notions of extremality related to conditions $\mathcal C$ and $\mathcal C_\nu$, and prove a relation between them. 

Here is a very general type of extremal solvability, which applies to interpolation data
\beq\label{verygen}
\la_j\in D \mapsto z_j\in E, \qquad j=1,\dots,n,
\eeq
where $D$ is a domain and $E$ is a connected subset of $\C^N$ for some $N$ (we have in mind sets $E$ that are either open or closed).
\begin{definition} \label{extlySolv}
The interpolation data \eqref{verygen} are {\em extremally solvable} if they are solvable but there do not exist an open neighbourhood $U$ of the closure of $D$ and a map $h\in\hol(U,E)$ such that the conditions 
\beq\label{interpGen1}
h(\la_j)=z_j \quad\mbox{ for } \quad  j=1,\dots,n,
\eeq
 hold.
\end{definition}
Thus a map $h\in\hol(\D,\Omega)$ is $n$-extremal if and only if, for any choice of $n$ distinct points $\la_1,\dots,\la_n\in\D$, the interpolation data $\la_j\in\D \mapsto h(\la_j)\in\Omega$ are extremally solvable.

Definition \ref{extlySolv} is natural from a geometric viewpoint, but it appears to be difficult to use in the context of Problem $I\Gamma$.   The following stronger notion has proved fruitful.
\begin{definition}\label{defCexly}
Let the interpolation data \eqref{Gdata} for Problem $I\Gamma$ satisfy condition $\mathcal C$.  The data {\em satisfy condition $\mathcal{C}$ extremally} (or {\em satisfy $\mathcal{C}_\nu$ extremally})
if there exists an $m\in\Schur$ (or  $m \in \B l_{\nu}$, respectively) such that the data 
\beq\label{active1}
\la_j \mapsto \Phi(m(\la_j), z_j), \quad j=1,\dots, n,
\eeq
are {\em extremally} solvable Nevanlinna-Pick data.   Alternatively, we say that the condition $\mathcal{C}(\la, z)$ (or $\mathcal{C}_\nu(\la, z)$)  {\em holds extremally}.   We say that $m\in\Schur$ or $\B l_\nu$ is an {\em auxiliary extremal} for the data \eqref{Gdata} if the data \eqref{active1} are extremally solvable.

We shall say that $\mathcal{C}_\nu$ is {\em active} or {\em holds actively and extremally} for the data \eqref{Gdata} if $\mathcal{C}_\nu(\la,z)$ holds extremally and there is a Blaschke product $m$ of degree exactly $\nu$ such that the data \eqref{active1} are extremally solvable.
\end{definition}
The conditions $\mathcal{C}_\nu$ were introduced in \cite{ALY12}.  In Definition \ref{defCexly} we do {\em not} assume that the interpolation data are solvable;  in fact one of the main questions that we confront is whether data that satisfy condition $\mathcal C$ extremally are necessarily solvable.
In the case that the data {\em are} solvable, however, we can ask how the extremal $\mathcal C$ condition relates to extremal solvability in the sense of Definition \ref{extlySolv}.

\begin{theorem}\label{stronger}
 Let   $\la_1, \dots,\la_n$  be distinct points in $\D$. If the  interpolation data
$$
\la_j \in\D \mapsto z_j  \in \G,  \quad j=1,\dots, n,
$$ 
are solvable and satisfy condition $\mathcal{C}$ extremally then the data are extremally solvable.
\end{theorem}
This will be proved in Section \ref{extsolvdata} (Theorem \ref{stronger4}), where we shall see further that other natural notions of extremal solvability are also weaker than the extremal $\mathcal C(\la,z)$ condition.

One could say that a function $h\in\hol(\D,\G)$ is {\em $3$-$\mathcal{C}$-extremal} if, for every triple $\la_1,\la_2,\la_3$ of distinct points in $\D$, the interpolation data $\la_j \mapsto h(\la_j)$ satisfies condition $\mathcal{C}$ extremally.  Theorem \ref{stronger} then shows that every $3$-$\mathcal{C}$-extremal map in $\hol(\D,\G)$ is $3$-extremal.
 The question as to how much stronger $3$-$\mathcal{C}$-extremality is than $3$-extremality relates to the conjecture \eqref{mainQ}: if it is true then the two notions coincide.  At present it remains open whether the conjecture \eqref{mainQ} is true.

\begin{remark}\label{Cn}
 If the $n$-point interpolation data  \eqref{Gdata} satisfy condition $\mathcal C$ extremally then they satisfy $\mathcal C_n$ extremally.  \rm

Trivially they satisfy $\mathcal C_n$.  Let $m\in\Schur$ be such that the Nevanlinna-Pick data \eqref{active1} are extremally solvable.  Since the $n$-point Nevanlinna-Pick data $\la_j \mapsto m(\la_j)$ are solvable, there exists a Blaschke product $\psi$ of degree at most $n$ such that $\psi(\la_j)=m(\la_j)$ for each $j$ (use induction and Schur reduction, or see \cite[Theorem 6.15]{AgMcC}).  Then
\[
\la_j \mapsto \Phi(\psi(\la_j),z_j), \quad j=1,\dots,n,
\]
are extremally solvable Nevanlinna-Pick data, and so the data $\la\mapsto z$ satisfy $\mathcal C_n$ extremally.
\end{remark}

When interpolation data \eqref{Gdata} for Problem $I\Gamma$ satisfy condition $\mathcal C_\nu$ extremally then, by Definition \eqref{Gdata}, they admit an auxiliary extremal $m\in\B l_\nu$.
It is far from the case that the auxiliary extremal $m$ is uniquely determined, or even that the degree $d(m)$ is unique for a particular set of data, as the following examples show. 

\begin{example}\label{Examples5_2} \rm  \cite[Examples 5.2]{ALY12} 
In each of these examples choose any three distinct points $\la_1,\la_2,\la_3\in\D$ and define $z_j$ to be $h(\la_j)$.\\
(1) Let $h(\la)=(2r\la,\la^2)$ where  $0 < r <1$.  Every degree 0 inner function $m \in\T$ is an auxiliary extremal for $\mathcal{C}_1$; there is no auxiliary extremal of degree $1$. Therefore in this case $\mathcal{C}_1$ holds extremally, but $\mathcal{C}_1$ is inactive.\\
(2)  Let 
\beq\label{alignex}
h(\la)= \left(2(1-r)\frac{\la^2}{1+r\la^3}, \frac{\la(\la^3+r)}{1+r\la^3}\right), \qquad \la\in\D.
\eeq
  The function $m(\la)=-\la$ is an auxiliary extremal for $\mathcal{C}_1$; there is no auxiliary extremal of degree $0$.  Here $q(\la)=-\la^2$. In this case $\mathcal{C}_1$ holds extremally and actively.\\
(3) Let $f$ be a Blaschke product of degree $1$ or $2$ and let $h=(2f,f^2)$.  Every  $m \in \B l_1$ is an auxiliary extremal and, for every $m$, we have $q=-f$.
\end{example}

\section{The main theorem} \label{maintheorem}

In this section we explain and motivate Theorem \ref{main}.   Let us recall the statement: \em

Let $\lambda_1,\lambda_2, \lambda_3$ be distinct points in $\mathbb{D}$
and let $z_1,z_2,z_3 \in \G$. 
The following statements are equivalent.
\begin{enumerate}
\item  There exists an aligned $\G$-inner function $h$ of degree at most $ 4$ such that
$h(\lambda_j) = z_j$ for  $ j=1, 2, 3$;

\item  condition $\mathcal{C}_1(\la, z)$ holds extremally and actively,
  and the associated Problem $\dia$ is solvable.
\end{enumerate}
\rm

We must explain the terms {\em aligned} and {\em Problem $\dia$}.  For the former
we need the notion of {\em royal nodes}.
\begin{definition}\label{defroyalnodes}
The {\em royal variety} $\mathcal {V} \subset \C^2$ is $\{(2z,z^2): z\in\C\}$.
A point $\la\in\C$ is a {\em royal node} of a rational $\G$-inner function $h$ if $h(\la) \in \mathcal V$.
\end{definition}
This is a specialization to rational functions of \cite[Definition 7.8]{ALY12}.  If $h=(s,p)$ then the royal nodes of $h$ are the solutions of the equation $s^2=4p$, and so there are $2d(p)$ of them, counting multiplicities and possible solutions at $\infty$.  Royal nodes lying in $\T$ are particularly important: if $\omega\in\T$ is a royal node of $h$ then $|s(\omega)|=2$ and furthermore the curve $h(\exp(it)), 0\leq t\leq 2\pi$, in $b\Ga$ touches the edge of the M\"obius band for $\exp(it)=\omega$; see Lemma \ref{signif} below.

It transpires that there are two qualitatively different types of degree $4$ $\G$-inner functions with $3$ or more royal nodes in $\T$, corresponding to different cyclic orderings of certain triples of points on the circle.  
\begin{definition} \label{defaligned}
Let $h=(s,p)$  be a rational $\G$-inner function.  We say that $h$ is {\em aligned} if  $h(\D)\subset\G$,  the degree of $h$ is at most $4$ and there exist at least $d(p)-1$ distinct royal nodes of $h$ in $\T$ and, if $d(p) =4$, there are distinct royal nodes $\omega_1,\omega_2,\omega_3$ of $h$ in $\T$ such that the points $\half s(\omega_1), \half s(\omega_2), \half s(\omega_3)\in\T$ are distinct and in the opposite cyclic order to $\omega_1,\omega_2,\omega_3$. 
\end{definition}
The significance of royal nodes in $\T$ is connected with cancellations in the function $\frac{2mp-s}{2-ms}$, as discussed in Section \ref{cancell}.

\begin{remark} \rm
The assumption that $h(\D)\subset\G$ is needed to exclude the possibility that $h$ map $\D$ into the topological boundary of $\G$ (that $h$ be `superficial' in the sense of \cite[Definition 8.1]{ALY12}).
\end{remark}

By \cite[Theorem 5.6]{AY04}, the  $2$-extremals in $ \hol(\D,\G)$ are aligned functions of degrees $1$ and $2$. 
 Example \ref{Examples5_2}(2) is an aligned function of degree $4$; some $\G$-inner functions of degree $4$ that are {\em not} aligned are given in Examples \ref{excaddy}.

In Section \ref{Eclasses} we give a characterisation of  aligned functions in terms  of the $\mathcal C_1$-extremality of $3$-point interpolation data generated by these functions. We show that  aligned functions  are $3$-extremals in $ \hol(\D,\G)$.

Given a $3$-point problem $I\Gamma$ that satisfies condition $\mathcal C_1$ extremally, in order to construct an interpolating function $h = (s, p)$ we aim first to find a unimodular rational function $p$ with suitable properties, then to define $s$ in terms of $p$, $m$ and $q$ and show that $(s, p)$ is the required interpolating function.  Clearly $p$ must be a solution of the Nevanlinna-Pick problem $\la_j\mapsto p_j, \, j=1,2,3$, but this is not enough -- it turns out that $p$ must also satisfy certain boundary interpolation conditions.  
The problem of finding a suitable $p$  comes down to the following.

\noindent {\bf Problem $\diamondsuit$} 
\em Given data $\la_j,s_j,p_j, \ j=1,2,3$, 
that satisfy condition $\mathcal{C}_1$  extremally with auxiliary extremal $m\in\aut\D$ 
find a Blaschke product $p$ of degree at most $4$ such that
\begin{equation}\label {def1p}
p(\la_j) = p_j, \qquad  j =1,2, 3,
\end{equation}
and
\begin{equation}\label{def2p}
p(\tau_\ell) = \overline m(\tau_\ell)^2, \qquad  \ell =1,\dots, d(mq),
\end{equation}
where the $\tau_\ell$ are the roots of the equation $mq(\tau)=1$ and $q$ is the unique function in the Schur class such that
\[
 q(\la_j) = \Phi(m(\la_j), s_j, p_j), \qquad j=1,2, 3.
\] 
\rm
The proof of Theorem \ref{main} is given in Sections \ref {Eclasses} to \ref{mainproof}.  

As is well known (for example, \cite{bgr,BD,bolot,geo,ChenHu}) there is a criterion of solvability of Nevanlinna-Pick problems like Problem $\dia$ in terms of the positivity and rank of a Pick matrix.   Combination of such a criterion with Theorem \ref{main} yields the following.

\begin{corollary}\label{mainbis}
Let $\lambda_1,\lambda_2, \lambda_3$ be distinct points in $\mathbb{D}$
and let $z_1,z_2,z_3 \in \G$. 
The following statements are equivalent.
\begin{enumerate}
\item  There exists an aligned $\G$-inner function $h$ of degree at most $ 4$ such that
$h(\lambda_j) = z_j$ for  $ j=1, 2, 3$;

\item  condition $\mathcal{C}_1(\la, z)$ holds extremally and actively,
  and if  $m\in\aut\D, \, q\in\B l_2$ and $\tau_1,\dots,\tau_{1+d(q)}$ are as in the statement of Problem $\dia$ then there exist positive numbers $\rho_1,\rho_2,\rho_3$ such that the $(4+d(q))$-square matrix $M=\bbm m_{ij} \ebm $ is positive semi-definite and of rank at most $4$, where 
\beq\label{bigPick}
m_{ij} = \left\{ \begin{array}{cll} \ds\frac{1-\bar p_ip_j}{1-\bar\la_i\la_j}  & \mbox{ if } & 1\leq i,j \leq 3, \\
		~&~  &~  \\
   \ds\frac{1-\bar p_i\bar m(\tau_{j-3})^2}{1-\bar\la_i\tau_{j-3}} & ~ & 1\leq i\leq 3, \, 4\leq j\leq 4+d(q),\\
		~&~  &~  \\
	\ds\frac{1-m(\tau_{i-3})^2 p_j}{1-\bar \tau_{i-3}\la_j} & ~ & 4\leq i\leq 4+d(q), \, 1\leq j\leq 3, \\
		~&~  &~  \\
	\ds\frac{1-m(\tau_{i-3})^2\bar m(\tau_{j-3})^2}{1-\bar\tau_{i-3}\tau_{j-3}} &~& 4\leq i,j\leq 4+d(q)\mbox{ and } i\neq j, \\
		~&~  &~  \\
	\ds\rho_{i-3} &~& 4\leq i=j\leq 4+d(q).
\end{array} \right .
\eeq
\end{enumerate}
\end{corollary}

\begin{remark} \rm
  In statement (1) of Theorem \ref{main} and Corollary \ref{mainbis} we suppose that there is an $m\in\mathcal B l_1$ of degree $1$ with certain properties.  The function $m$ is not uniquely determined; each choice of $m$ generates a different Problem $\dia$.  A consequence of the theorem is that if there is some $m$ for which Problem $\dia$ has a solution, then the same holds for {\em all} auxiliary extremals $m$ of degree $1$.  See Example \ref{Examples5_2}(3) above for illustration.
\end{remark}

\section{The classes $\mathcal E_{\nu n}$ of rational functions}\label{Eclasses}

In this section we recall some results from \cite{ALY12} about some classes of rational $\G$-inner functions.
The following statement is proved in \cite[Proposition 5.1]{ALY12}; it follows easily from the properties of solutions of extremally solvable Nevanlinna-Pick problems.
\begin{proposition} \label{equivCnu}
For any $\Gamma$-interpolation data $\la_j \mapsto z_j, \ j=1,\dots, n$ and $\nu \ge 0$, the following conditions are equivalent.

{\rm (i)} $\mathcal{C}_\nu(\la, z)$  holds extremally;

{\rm (ii)}  $\mathcal{C}_\nu(\la, z)$  holds and there exist
$ m \in \B l_\nu$ and $q \in \B l_{n -1}$ such that 
\begin{equation}
\label{Phiq}
\Phi(m(\la_j), z_j)= q(\la_j), \quad j =1, \dots, n.
\end{equation}

Moreover, when condition {\rm (ii)} is satisfied for some $m \in \B l_\nu$, there is a {\em unique}  $q \in \B l_{n -1}$ such that equations \eqref{Phiq} hold. If, furthermore, the 
$\Gamma$-interpolation data $\la_j \mapsto z_j, \ j=1,\dots, n$, are solvable by an analytic function $h:\D\to\Gamma$ then
\beq\label{propq}
\Phi\circ (m,h) = q.
\eeq
\end{proposition}
Thus, when (ii) holds, if $h=(s,p)$ then
\beq\label{q=2mpetc}
\frac{2mp - s}{2 - ms} = q.
\eeq

Proposition \ref{equivCnu} leads us to consider some classes $\E_{\nu n} \subset \hol(\D,\Gamma)$  of rational functions. 
\begin{definition}\label{defEclass} For  $\nu \ge 0$, $n \ge 1$, we say that a  function $h=(s,p) $ is in  $\E_{\nu n}$ (or in $\widetilde{\E_{\nu n}}$) if $h\in\hol(\D,\Gamma)$ is rational and there exists $m \in \B l_{\nu}$
(or a Blaschke product  $m$ of degree $\nu$, respectively) such that 
$$
\frac{2 m p -s}{2- m s} \in \B l_{n-1}.
$$
\end{definition}
Thus Proposition \ref{equivCnu} states that if $h$ is a solution of  Problem $I\Gamma$ with data that satisfies condition $\mathcal C_\nu$ extremally, then $h\in \E_{\nu n}$.  A function in $ \E_{\nu n}$ is rational of degree at most $2n-2$ and is necessarily { $\Gamma$-inner} \cite[Theorem 7.3]{ALY12}.
\begin{proposition} For $n \ge 1$  and $\nu \ge n$ 
$$
 \E_{\nu n}= \E_{n n}.
$$
\end{proposition}
\begin{proof} Let  $h=(s,p) $ be in  $\E_{\nu n}$. There exists $m \in \B l_{\nu}$ such that 
$$
\frac{2 m p -s}{2- m s} \in \B l_{n-1}.
$$
Let $\la_1, \dots, \la_n$ be distinct points in $\D$. Then 
the $n$-point interpolation data
\beq\label{Gdata_n}
\la_j \mapsto \frac{2 m(\la_j) p(\la_j) -s(\la_j)}{2- m(\la_j) s(\la_j)}, \quad j=1,\dots, n,
\eeq
are  extremally solvable Nevanlinna-Pick data.  
 Since the $n$-point Nevanlinna-Pick data $\la_j \mapsto m(\la_j)$ are solvable, there exists a Blaschke product $\psi$ of degree at most $n$ such that $\psi(\la_j)=m(\la_j)$ for each $j$ (use induction and Schur reduction, or see \cite[Theorem 6.15]{AgMcC}).  Then
\[
\la_j \mapsto \Phi(\psi(\la_j),s(\la_j), p(\la_j)), \quad j=1,\dots,n,
\]
are extremally solvable Nevanlinna-Pick data.

By Proposition \ref{equivCnu},
\[
\frac{2\psi p - s}{2 - \psi s}=\frac{2mp - s}{2 - ms}  \in \B l_{n-1}.
\]
Since $\psi \in \B l_{n}$, the function $h=(s,p) $ is in  $\E_{n n}$.
\end{proof}

Consider a function $h=(s,p)\in\widetilde {\mathcal{E}_{13}}$ of degree $4$.  Let $m\in\aut\D$ be an auxiliary extremal, so that equation \eqref{q=2mpetc} holds and $q$ has degree at most $2$.  Since $d(mp)=5$ there must be at least $3$ cancellations between the numerator and denominator in equation \eqref{q=2mpetc}.  An understanding of these cancellations will be important in the sequel;  we recall some results from \cite{ALY12} about them.

\begin{lemma}\label{signif}
If $(s,p)$ is a non-constant rational $\G$-inner function then
\begin{enumerate}
\item  the points in $\Delta$ at which $|s|=2$ are precisely the royal nodes of $(s,p)$ in $\T$, and
\item for any finite Blaschke product $m$,  the rational function $\frac{2mp-s}{2-ms}$ has a cancellation at $\zeta\in\C$ if and only if $\zeta$ is a royal node of $(s,p)$ in $\T$ and $m(\zeta) = \half\overline{s(\zeta)}$.  Moreover, when this is so there is exactly one cancellation at $\zeta$.
\end{enumerate}
\end{lemma}
The first assertion is Lemma 7.10, the second is Theorem 7.12 in \cite{ALY12}.

\begin{proposition}\label{aligned_nodes}
A function  $h \in \hol(\D,\G)$
 is aligned if and only if  $h\in \widetilde{\mathcal E_{13}}$.
\end{proposition}

\begin{proof} Suppose $h=(s,p)\in \widetilde{\mathcal E_{13}}$.  By \cite[Theorem 7.3]{ALY12} any function in ${\mathcal E_{\nu n}}$ is  rational, $\G$-inner  and has degree at most $2n-2$.  Thus $d(p) \leq 4$.
By \cite[Corollary 6.10]{ALY12},  $s$ is rational of degree at most $4$ and has the same denominator as $p$.  Let 
\beq\label{qform}
q=\frac{2mp-s}{2-ms};
\eeq
then, by assumption, $d(q)\leq 2$.   If $d(p)=4$ then, since  $d(m)=1$, the degree of $q$ is $5 $ minus the number of cancellations in the right hand side of equation \eqref{qform}, and therefore there are at least $3$ cancellations.  By Lemma \ref{signif}, there is exactly one cancellation at each royal node $\omega$ of $h$ at which $m(\omega)= \half \overline{s(\omega)}$.  Hence there are $3$ distinct royal nodes $\omega_j$ of $h$ such that the $3$ points $\half\overline{s(\omega_j)}=m(\omega_j)$ are distinct and are in the same cyclic order as the $\omega_j$.  Thus $h$ is aligned.

If $d(p) =3$ then similar reasoning shows that there are at least $2$ cancellations in the right hand side of equation \eqref{qform} and hence $h$ has two royal nodes in $\T$.  Likewise if $d(p)=2$ then $h$ has a royal node in $\T$.  In all cases $h$ is aligned.

Conversely, suppose that $h$ is aligned; then $d(p) \leq 4$.  We prove that $h\in \widetilde{\mathcal E_{13}}$ in the case that $d(p)=4$.  By hypothesis there are distinct royal nodes $\omega_1,\omega_2,\omega_3$ of $h$ in $\T$ such that the points $\half s(\omega_1), \half s(\omega_2), \half s(\omega_3)\in\T$ are distinct and in the opposite cyclic order to $\omega_1,\omega_2,\omega_3$.   It follows that there exists $m\in\aut\D$ such that $m(\omega_j)= \half\overline{s(\omega_j)}$ for $j=1,2,3$.  By Lemma \ref{signif}, there are $3$ cancellations in the right hand side of equation \eqref{qform}, and hence $q\in\B l_2$.  Thus $h\in\widetilde{\mathcal E_{13}}$.  In the case that $d(p)<4$ the function $h$ has at least $d(p)-1 \leq 2 $ royal nodes in $\T$ and so there exists $m\in\aut\D$ such that $m(\omega)=\half \overline{s(\omega)}$ at each of these nodes $\omega$.  Then there are $\max\{0,d(p)-1\}$ cancellations in the right hand side of equation \eqref{qform}, and hence $q\in\B l_2$.  Thus $h\in\widetilde{\mathcal E_{13}}$.  
\end{proof}
\begin{corollary}\label{aligned_extremal}
Every aligned  function is $3$-extremal in  $ \hol(\D,\G)$.
\end{corollary}
\begin{proof}
By \cite[Theorem 9.1]{ALY12}, any function in $\E_{1 3}$ which is not  superficial  is $3$-extremal. Recall that a superficial function maps into the topological boundary of $\Gamma$. By Definition \ref{defaligned}, an aligned function $h$ belongs to $ \hol(\D,\G)$ and so is not superficial. By Proposition \ref{aligned_nodes}, $h\in \widetilde{\mathcal E_{13}}$. Thus $h$ is  $3$-extremal in $ \hol(\D,\G)$.
\end{proof}

The following preparatory results will be used in the proof of Theorem \ref{main}.
There is a special case of Problem $I\Gamma$ in which condition $\mathcal{C}_0$ {\em is} sufficient, for any number of nodes.
\begin{proposition}
\label{flat}
Let $\lambda_1, \dots, \lambda_n$ be distinct points in $
\mathbb{D}$ and let $z_j=(s_j, p_j) \in \G, \   j =1,\dots, n$. 
If condition $\mathcal{C}_0(\la,z)$ holds and
the Nevanlinna-Pick
problem with data $\lambda_j \mapsto p_j$ is {\em extremally}
solvable then
$$
\la_j \mapsto z_j, \quad 1 \le j \le n,
$$
are solvable $\Gamma$-interpolation data.
\end{proposition}
This result is \cite[Theorem 5.2]{AY04}.

 Pick's Theorem enables us to recast the  condition
$\mathcal{C}_\nu$ 
as the positivity of a pencil of matrices.
\begin{proposition}
\label{uppick}
If
$$
\la_j \mapsto z_j=(s_j,p_j), \quad 1 \le j \le n,
$$
are interpolation data for $\Gamma$ then condition 
$\mathcal{C}_\nu(\la_1,\dots,\la_n,z_1,\dots,z_n)$ holds
if and only if,
for every Blaschke product $\up$ of degree at most $\nu$,
{\small
\begin{equation}
\label{pick}
\left[\frac{1-\up(\lambda_i) p_i \bar p_j \overline
\up(\lambda_j) - \frac{1}{2} \up(\lambda_i) (s_i - p_i \bar s_j) -
\frac{1}{2} (\bar s_j - \bar p_j s_i) \overline \up
(\lambda_j)
  -\tfrac 14   (1-\ups(\la_i) \bar\ups(\la_j) ) s_i \bar s_j
}{1-\lambda_i \overline \lambda_j}\right]^n_{i, j=1}
\ge 0.
\end{equation}
}
\end{proposition}
The details are given in \cite[Theorem 4.5]{ALY12}.

\section{Calculation of interpolating functions} \label{calculation}
In this section we present two lemmas on the construction of interpolating functions for data (with any number of nodes) that satisfy condition $\mathcal C$ extremally.
Suppose that $n \ge 3$ and data $\la_j \mapsto (s_j, p_j)$, $ j=1,\dots, n$, satisfy $\mathcal{C}_{n-2}$ extremally, with auxiliary extremal $m \in \B l_{n-2}$.
Here is a strategy for constructing an interpolating function $h = (s, p)$:

(1) find a unimodular rational function $p$ with suitable properties, 

(2)  define $s$ in terms of $p, m$ and $q$ as in Proposition \ref{equivCnu}, so that equation (\ref{propq}) holds, and then 

(3) show that $(s, p)$ is the required interpolating function.

A delicate part of the process is to show that $(s, p)(\D)\subset \G$.  The next result describes the construction of the functions $s$ and $p$ and shows that if $|s|\leq 2$ on $\T$ then it will follow that $(s,p)(\D)\subset\G$.

 For a function $f$ on a subset of the complex plane $\mathbb{C}$
 write $\bar f(z) = \overline{f(z)}.$

\begin{lemma} \label{construct}
Suppose that data 
\[
\la_j \mapsto z_j=(s_j, p_j),\; \;\;  j=1,\dots, n,
\]
satisfy $\mathcal{C}_{n-2}$ extremally, with auxiliary extremal $m \in \B l_{n-2}$ and let $q$ be the unique member  of $\B l_{n-1}$ such that 
\begin{equation}
\label{Phiq-const}
\Phi(m(\la_j), z_j)= q(\la_j), \quad j =1, \dots, n.
\end{equation}
Suppose that $p$ is a rational inner function such that
\begin{equation}\label{p-at-la}
p(\la_j) = p_j, \qquad  j=1,\dots,  n,
\end{equation}
and that the function $s$ is defined by the equation
\begin{equation}\label{defs}
s = \frac{2(mp - q)}{1-mq}.
\end{equation}
If $s$ satisfies $|s| \le 2$ on $\T$ then the function $h = (s,p)$ is an analytic
function from $\D$ to $\Gamma$ such that
$$
h(\la_j) = (s_j, p_j), \qquad  j=1,\dots,  n.
$$
\end{lemma}
\begin{proof}
Note that $h: \mathbb{D} \to \mathbb{C}^2$ is analytic and
\begin{equation} \label{sj}
s(\lambda_j) = 2 \ \frac{m(\lambda_j) p_j - q(\lambda_j)}{1 -
m(\lambda_j) q(\lambda_j)} = s_j
\end{equation}
(by choice of $q$ satisfying (\ref{Phiq-const})).  
Thus
$h(\lambda_j) = z_j$ for $j=1,\dots, n$. We must prove that $h(\D) \subset \Gamma$.

  For $\lambda
\in \mathbb{T}$ we have $|p(\lambda)| = 1 =
|q(\lambda)|=|m(\lambda)|$ and
\begin{eqnarray}
\label{hGinner}
\bar s p(\lambda) &=& 2 \ \frac{\overline m \bar p(\lambda) -
\bar q (\lambda)}{1 - \overline m \bar q(\lambda)} p(\lambda)\\
&=& 2 \  \frac{\overline m(\lambda) - p\bar
q (\lambda)}{1 - \overline m \bar q(\lambda)} = 2 \
   \frac{q(\lambda) - mp (\lambda)}{ m q(\lambda)-1} \nonumber
\\ &=& s(\lambda). \nonumber
\end{eqnarray}
Hence, by \cite[Proposition 3.2]{ALY12}  and the fact that
  $|s(\lambda)| \le 2$ for all $\lambda \in
\mathbb{T}$, we have
  $h(\lambda) \in \Gamma$ for all $\lambda \in
\mathbb{T}$.  Thus, for
fixed $z\in \mathbb{D}$, the analytic mapping $\lambda \mapsto
\Phi(z, h(\lambda))$ is bounded by 1 on $\mathbb{T}$ and so,
by the Maximum Principle, is bounded by 1 on $\mathbb{D}$.  By \cite[Proposition 3.2]{ALY12}, $h(\lambda) \in \Gamma$ for all $\lambda \in
\mathbb{D}$. \hfill\qed
\end{proof}

The boundedness of $s$ on $\T$ places a restriction on $p$.

\begin{lemma}\label{pattau}
Let $m$, $q$ and $p$ be rational inner functions and let $s$ be defined by the equation 
\begin{equation}\label{defs-const}
s = \frac{2(mp - q)}{1-mq}.
\end{equation}
Let the solutions of the equation $mq(\tau) =1$ be the points
\begin{equation}\label{deftau}
\tau_\ell,\qquad  \ell=1,\dots, d(mq).
\end{equation}
If $s$ is bounded on $\T$ then
\begin{equation}\label{p-at-tau}
p(\tau_{\ell}) = \overline{m(\tau_{\ell})}^2, \qquad \ell =1,\dots, d(mq).
\end{equation}
\end{lemma}
\begin{proof}
If $\tau$ is a solution of $mq(\tau) =1$, the boundedness of $s$ implies that the numerator $mp -q$ on the right hand side of equation (\ref{defs-const}) vanishes at $\tau$.
  Note that since $mq$ is a Blaschke product any solution of the
equation  $mq(\tau) =1$ lies on $\mathbb{T}$. Hence $q(\tau)=\overline{m(\tau)} $. Thus
\[
(mp -q)(\tau)= m(\tau) p(\tau) -\overline{m(\tau)} =0,
\]
which implies that 
$p(\tau) = \overline{m(\tau)}^2$.
\hfill\qed
\end{proof}

Proposition \ref{equivCnu} and Lemmas \ref{construct} to \ref{pattau} show that if our construction strategy is to succeed the first step must be to find a rational inner function $p$ satisfying interpolation conditions \eqref{p-at-la} at the $\la_j \in \D$ and \eqref{p-at-tau} at the $\tau_\ell\in\T$.  Consideration of the number of degrees of freedom suggests we should seek $p$ that is a Blaschke product of degree at most $2n - 2$.

\section{Properties of interpolating functions} \label{property}

To establish a sufficient condition for three-point interpolation in
$\hol(\D, \Gamma)$ we shall need  some technical observations.

A rational function $f$ is {\em unimodular} if $|f(z)| =1$ for all $z \in \T$.

\begin{lemma}\label{psiinner}
Let
$$
\lambda_i \mapsto z_i, \qquad  i=1,\dots,  n,
$$
be solvable Nevanlinna-Pick data.  If $\psi$ is a rational
function of degree $n$ which is unimodular on $\mathbb{T}$ and
satisfies $\psi(\lambda_i) = z_i , \ i=1,\dots,  n$,  then $\psi$
is a Blaschke product.
\end{lemma}

\begin{proof}
Consider the case $n=1$.  We must have $\psi = cB_\alpha$ or
$\psi=c/B_\alpha$ for some $c\in \mathbb{T}$ and $\alpha \in
\mathbb{D}$.  In the latter case we have
$$c = B_\alpha \psi(\lambda_1) = B_\alpha(\lambda_1) z_1$$
and so $|c| <1$, which is a contradiction.  Hence $\psi = cB_\alpha$,
a Blaschke product.

Now consider the case $n>1$ and suppose the result known for
$n-1$.  Let $\psi_1 = B_{z_1} \circ \psi/ B_{\lambda_1}$.  Then
$\psi_1$ is rational of degree $n-1$ and unimodular on
$\mathbb{T}$.  Furthermore,
$$
\psi_1(\lambda_i) =
\frac{B_{z_1}(z_i)}{B_{\lambda_1}(\lambda_i)},\qquad i= 2,\dots,n.
$$
By the standard Schur reduction process \cite{W},
$$
\lambda_i\mapsto
\frac{B_{z_1}(z_i)}{B_{\lambda_1}(\lambda_i)},\qquad  i=2,\dots, n,
$$
are solvable Nevanlinna-Pick data.  By the inductive hypothesis,
$\psi_1$ is a Blaschke product.  Thus
$$
\psi = B_{-z_1} \circ (\psi_1 B_{\lambda_1})
$$
is also a Blaschke product.  Hence the assertion holds for all $n\in
\mathbb{N}$.
\end{proof}
\begin{lemma}
\label {propsp}
Let $p$ be a rational function, let $m,\ q$ be Blaschke products and let $s$ be defined by equation {\rm(\ref{defs})}.  Then
\beq \label {royal}
s^2 - 4p = \frac{4(m^2p-1)(p-q^2)}{(1-mq)^2}.
\end{equation}
Let
$$
p = \frac{N_p}{D_p}
$$
where $N_p, \ D_p$ are coprime polynomials.  If further
\begin{eqnarray}
\label {sometau}
p(\tau_\ell) &=& \bar m(\tau_\ell)^2,  \quad  \ell=1,\dots,  J,
\end{eqnarray}
for some $J \le d(mq)$, where the $\tau_\ell$ are distinct solutions of $mq=1$,
  then, for some polynomial $R$ of degree at most $d(p)+d(mq)-J$,
\begin{equation}
\label{exprs}
s=\frac{R}{D_p \Pi}
\end{equation}
where
\begin{equation} \label{defPi}
\Pi(\la)= \prod_{J < \ell \le d(mq)} \la - \tau_\ell.
\end{equation}
\end{lemma}
\begin{proof}
The identity (\ref{royal}) is straightforward.   In a self-explanatory notation, we have
$$
s = \frac{2}{D_p}\ \frac{N_m N_p D_q - N_q D_m D_p}{D_m D_q -
N_m N_q}
$$
as a ratio of polynomials.  By choice of $p$ to satisfy equation
(\ref{sometau}), the factors $\la - \tau_\ell,\  \ell=1,\dots, J,$
cancel on the right hand side, and we obtain the expression (\ref{exprs}). 
\end{proof}

The example below shows that the rational functions $s$ and $p$ from Lemma \ref{propsp} may have different degrees.

\begin{example} \label{surprise}\rm Let $a \in \D \setminus  \{ 0\}$ and let
\[
h(\la) = (s(\la),p(\la)) = \left(\frac{c \la}{1-\bar{a} \la}, \frac{\la(\la-a)}{1-\bar{a} \la} \right)
\]
for some $c \in \R$ such that $|c| \le 2(1- |a|)$.
Clearly $h\in\hol(\D,\Gamma)$ and $d(s) =1$ and $d(p)=2$; see
 \cite[Example 6.7]{ALY12}.
\end{example}

\section{Cancellations in some rational functions} \label{cancell}

Underlying the technical results in this paper is a study of cancellations in certain rational functions.  Corresponding to a function $h$ which is a candidate for a solution to an interpolation problem we introduce
 the function $\ph_\up=\Phi \circ (\up, h)$, where  $\upsilon$  is  Blaschke factor.  An understanding of  cancellations in $\ph_\up$ will enable us to show that $h$ is analytic in $\D$.  We have previously studied such
cancellations in  \cite[Section 7.2]{ALY12}.

\begin{lemma}\label {getphiu}
Let $m, \ p, \ q$ and $ \upsilon$ be unimodular rational functions. The function
  $\varphi_\ups$ defined by
\begin{equation}\label{phiu}
\varphi_\up =  \Phi\circ(\up,s,p)=\frac{2\up p -s}{2 - \up s},
\end{equation}
where 
\begin{equation}\label{defs2}
s = \frac{2(mp - q)}{1-mq},
\end{equation}
is a unimodular rational function.  
\end{lemma}

\begin{proof}
At any point of $\T$ we have $|\up p|=1$ and (see equations (\ref{hGinner}))
$s = \bar s p$, and so
$$
|\varphi_\up| = \left|\frac{2 \up p -s}{2 - \up s}\right| = \left|
\frac{2-\overline \up \bar p s}{2 - \up s}\right| = \left| \frac{2 -
\overline \up \bar s}{2 - \up s }\right| =1. 
$$
\end{proof}

We plan to use Lemma \ref{psiinner} to show that, for suitable $\up$, $\ph_\up$ is analytic in $\D$, that is, it is a Blaschke product. We therefore need to know the degree of $\ph_\up$, and to this end  we need to study cancellations.  We shall say that  $\ph_\up$ {\em has $n$ cancellations} at a point $\al$  if both numerator and denominator in the right hand side of equation (\ref{phiu}) vanish at $\al$ to order $n$. We summarise the possibilities for cancellations.

\begin{lemma}\label{cancels}
Let $m, \ p, \ \upsilon$ and $q$ be unimodular rational functions.   Let $s, \ \ph_\up$ be given by equations $(\ref{defs2}), (\ref{phiu})$ respectively.
\begin{enumerate}
\item [\rm ($\al$)] If  $\ph_\up$ has a cancellation at $\al$ then either $m^2p(\al)=1$ or $p(\al)=q(\al)^2$.
\item [\rm ($\si$)] Let $m^2p(\si)=1$ and $mq(\si)\neq 1$:  then $\ph_\up$ has a cancellation at $\si$ if and only if $\up(\si)=m(\si)$.  If $\si$ is a zero of $m^2p-1$ of order $\nu$ and is a zero of $\up -m$ of order $n\leq \nu$ then $\ph_\up$ has $n$ cancellations at $\si$.
\item [\rm ($\beta$)] Let $mq(\beta)\neq 1$ and  $p(\beta)=q(\beta)^2\neq 0$: then $\ph_\up$ has a cancellation at $\beta$ if and only if $\up(\beta)= -1/q(\beta)$.
\item [\rm ($\beta\si$)]  Let $mq(\beta)\neq 1, \ m^2p(\beta)=1$ and $p(\beta)=q(\beta)^2$: then $\ph_\up$ has a double cancellation at $\beta$ if and only if $\up(\beta)= m(\beta)$ and $\up'(\beta) = -\half m(\beta)^3p'(\beta)$.
\item[\rm($\beta\beta$)] Let $mq(\beta)\neq 1, q(\beta)\neq0$ and let $\beta$ be a double zero of $p-q^2$.  Then $\ph_\up$ has a double cancellation at $\beta$ if and only if $\up(\beta)= -1/q(\beta)$ and $\up'(\beta) = q'(\beta)/q(\beta)^2$.
\item [\rm ($\tau$)] Let  $\tau$ be a simple zero of $mq-1$ and a double zero of $p-q^2$: then $\ph_\up$ has a cancellation at $\tau$ if and only if $\up(\tau)=-m(\tau)$.
\item[\rm($\tau\beta$)] Let  $\tau$ be a simple zero of $mq-1$ and a triple zero of $p-q^2$: if $q'(\tau)\neq 0$ then $\ph_\up$ does not have a double cancellation at $\tau.$
\end{enumerate}
\end{lemma}
\begin{proof}
($\al$) Suppose there is a cancellation at $\al$, that is, $(2\up p -s)(\al)=0=(2-\up s)(\al)$.
Then 
\[
 (s^2-4p)(\al)= (s^2-2\up s p)(\al) = -s(\al)(2\up p -s)(\al)=0.
\]
It follows from equation (\ref{royal}) that either $m^2p(\al)=1$ or $p(\al)=q(\al)^2$.

\noindent ($\si$)  Note the identities
\begin{eqnarray*}
2mp -s &=& -\frac{2q(m^2p - 1)}{1-mq},\\
2-ms &=&  -2\frac{m^2p-1}{1-mq}.
\end{eqnarray*}
It follows that
\begin{eqnarray*}
2\up p -s &=& -\frac{2q(m^2p - 1)}{1-mq} + 2p(\up -m),\\
2-\up s &=&  -2\frac{m^2p-1}{1-mq} - s(\up - m).
\end{eqnarray*}

By assumption, $m^2p(\si)=1$ and $mq(\si)\neq 1$, thus, if $2\up p-s$ and $2-\up s$ vanish at $\si$ then so does $\up -m$ (note that $p(\si)\neq 0$).  Conversely, if  $\up - m$ vanishes at $\si$
with multiplicity $n \le \nu$, then so do $p(\up - m)$
and  $s(\up - m)$, and the numerator and denominator
of $\varphi_\up$ have $n$ cancellations at $\si$.

\noindent ($\beta$) is easy.

\noindent ($\beta\si$)  The assumptions imply that $mq(\beta)=-1$.  On differentiating equation (\ref{defs}) we find
\beq\label{derivs}
s'= 2\frac{m'(p-q^2) + (m^2p-1)q' + mp'(1-mq)}{(1-mq)^2}.
\eeq
Hence
\[
s'(\beta)= 2\frac{mp'}{1-mq}(\beta)=(mp')(\beta).
\]
It follows that $2\up p - s= 2-\up s = (2\up p - s)'=(2-\up s)'=0$ at $\beta$ if and only if $\up(\beta)=m(\beta)$ and $\up'(\beta)=-\half (m^3p')(\beta)$.

\noindent ($\beta\beta$)  By ($\beta$) above, there is one cancellation at $\beta$ if and only if $\up(\beta)=-1/q(\beta)$.    We have $p(\beta)=q(\beta)^2, p'(\beta)=2qq'(\beta)$, from which it follows that $s(\beta)=-2q(\beta), \, s'(\beta)=-2q'(\beta)$.  From these equations it is straightforward to calculate that $(2\up p-s)'(\beta)=(2-\up s)'(\beta)=0$ if and only if
$\up'(\beta)=q'(\beta)/q(\beta)^2$.

\noindent ($\tau$)   We have $m(\tau) \neq 0$ and $p'(\tau)=2qq'(\tau) = 2q'(\tau)/m(\tau)$.  By L'H\^{o}pital's Rule 
 \[
s(\tau)=-2\frac{m'p+mp'-q'}{m'q+mq'}(\tau) =-2\frac{m'p+q'}{m'mp+mq'}(\tau)= -2/m(\tau),
\]
 and the assertion follows easily. 

\noindent ($\tau\beta$)  Suppose that $\ph_\up$ does have a double cancellation at $\tau$.  We have $\up(\tau)=-m(\tau) \neq 0, \, p'(\tau)=2q'(\tau)/m(\tau)$ and
\[
p''(\tau)=(2qq')'(\tau)= 2q'(\tau)^2+2q''(\tau)/m(\tau).
\]
Since $\tau$ is a triple zero of $p-q^2$,
\[
\frac{p-q^2}{(1-mq)^2}(\la) \to 0 \qquad \mbox{as  } \la \to \tau,
\]
and so, by equation (\ref{derivs}),
\[
s'(\tau)= \lim_{\la\to\tau}\frac{fq'+mp'}{1-mq}(\la)
\]
where
\[
f=\frac{m^2p-1}{1-mq} = m^2\frac{p-q^2}{1-mq}-1-mq.
\]
As $\la\to\tau$ we have
\[
f(\la) = -2 - (mq)'(\tau)(\la-\tau)+O(\la-\tau)^2.
\]
Thus $f(\tau)=-2, \, f'(\tau)=-(mq)'(\tau)$ and
\begin{eqnarray*}
s'(\tau) &=& \frac{(fq'+mp')'}{-(mq)'}(\tau) = \frac {(mq)'q'+2q''-m'p'-mp''}{(mq)'}(\tau) \\
          &=& \frac{m'q'/m+m(q')^2+2q''-2m'q'/m-2m(q')^2-2q''}{m'/m+mq'}(\tau)\\
             & =& \frac{-m'q'-m^2(q')^2}{m'+m^2q'}(\tau) = -q'(\tau).
\end{eqnarray*}
Thus
\[
(2\up p -s)'(\tau) = 2\frac{\up'}{m^2}(\tau) +2(-m)2\frac{q'}{m}(\tau)+q'(\tau)
       = \frac{2\up'-3m^2q'}{m^2}(\tau).
\]
Hence $(2\up p-s)'(\tau)=0$ if and only if $\up'(\tau)= \tfrac 32 (m^2q')(\tau)$.

On the other hand, since $s(\tau)=-2/m(\tau),$
\[
(2-\up s)'(\tau)= -(\up' s + \up s')(\tau)=\frac{2\up' -m^2q'}{m}(\tau),
\]
so that $(2-\up s)'(\tau)=0$ if and only if $\up'(\tau)= \tfrac 12 (m^2q')(\tau)$.  Thus $(2\up p-s)'$ and $(2-\up s)'$ cannot simultaneously vanish at $\tau$ when $q'(\tau)\neq 0$.
\hfill $\Box$
\end{proof}

We can sharpen the results of Lemma \ref{cancels} with the aid of the following elementary notion.
\begin{definition}\label{phasar-derivative}
 For any differentiable function $f: \T \to \C \setminus \{0\}$
 the {\em phasar derivative} of $f$ at $z= \e^{\ii \ta} \in \T$
is  the derivative with respect to $\ta$ of the argument of $f(\e^{\ii \ta})$ at $\ta$; we denote it by $Af(z)$.
\end{definition}
Thus, if $f(\e^{\ii \ta}) = R(\ta) \e^{\ii g(\ta)}$ where $g(\ta) \in\R$ and $R(\ta) >0$ then $g$ is differentiable on $[0,2\pi)$ and the phasar derivative of $f$ at $z= \e^{\ii \ta} \in \T$ is equal to
\begin{equation}\label{phasar-deriv-arg}
Af(\e^{\ii \ta}) = \frac{d}{d \ta} \arg f(\e^{\ii \ta}) = g'(\ta).
\end{equation}
Clearly, for differentiable functions $\psi, \phi : \T \to \C \setminus \{0\}$ and
for any $c \in \C \setminus \{0\}$, we have
\beq\label {phasar-deriv-sum}
A(\psi \phi ) =A \psi +A \phi \; \quad \text{and} \; \quad A(c \psi  ) =A \psi.
\eeq

The result below on properties of phasar derivatives are simple; they can be found in \cite[Section 7.1]{ALY12}.
\begin{proposition}\label{AofInner}  Let $\phi: \T \to \C \setminus \{0\}$ be a rational inner function. Then, for all $\lambda \in \T$,
\begin{equation} \label{phasar-deriv-inner}
 A\phi(\lambda)= \lambda \frac{\phi'(\lambda)}{\phi(\lambda)} \;\;\text{and}\; \; A \phi (\lambda) >0.
\end{equation}
\end{proposition} 

In the applications of Lemma \ref{cancels} below, $\up, q$ are Blaschke products.  In such cases there are few possibilities for double cancellations on the unit circle.
\begin{corollary}\label{fewdouble}
Let $\up$  be a finite Blaschke product and let $m, q$ and $p$ be unimodular rational functions.  Let $s, \ph_\up$ be as in Lemma {\rm \ref{cancels}}.
\begin{enumerate}
\item[\rm($\beta\si$)]  Let $mq(\beta)\neq 1, \ m^2p(\beta)=1$ and $p(\beta)=q(\beta)^2$ for some $\beta\in\T$.  Then $\ph_\up$ has a double cancellation at $\beta$ if and only if $\up(\beta)= m(\beta)$ and $A\up(\beta) = -\half Ap(\beta)$.  In particular, $\ph_\up$ has no double cancellation at $\beta$ if $p$ is inner.
\item[\rm($\beta\beta$)] Let $mq(\beta)\neq 1, q(\beta)\neq0$ and let $\beta$ be a double zero of $p-q^2$ for some $\beta\in\T$.  Then $\ph_\up$ has a double cancellation at $\beta$ if and only if $\up(\beta)= -1/q(\beta)$ and $A\up(\beta) = -Aq(\beta)$.  In particular, $\ph_\up$ has no double cancellation at $\beta$ if $q$ is inner.
\end{enumerate}
\end{corollary}
\begin{proof}
The conditions on $A\up(\beta)$ are simply restatements of the corresponding items in Lemma \ref{cancels} in terms of phasar derivatives.  The impossibility of double cancellations in the case of inner $p, q$ follows from the fact that $Af > 0$ on $\T$ for any finite Blaschke product $f$. \hfill $\Box$
\end{proof}

\section{Snares} \label{snlemma}
To establish a sufficient condition for three-point interpolation
$\D \to \Gamma$ we shall need a topological lemma on multifunctions
in order to prove some delicate boundedness properties.

If $X$ and $Y$ are topological spaces then a {\em multifunction} from $X$ to $Y$ is defined to be a mapping from $X$ to the set of all subsets of $Y$.
Such a multifunction $S$ is said to be {\em upper
semi-continuous} if $\{\la: {S}(\la) \subset U \}$ is open in $X$ for every
open set $U$ in $Y$.  If $S_1$ and $S_2$ are multifunctions 
from $X$ to $Y$ then so is $S_1 \cup S_2$, where 
$(S_1 \cup S_2)(\la) \df S_1(\la) \cup S_2(\la)$ for $\la\in X$. 
\begin{definition}
\label{defsnare}
A {\em snare} is  a multifunction $ S$ from a subset $X$ of $\D$ to $\C^* \setminus 2\D$ with the following properties:
\begin{enumerate}
\item[\rm (1)] $X$ is a connected open subset of $\D$ and the closure of $X$ in
$\C$ contains $\T$;
\item [\rm (2)] $ S(\la)$ is a compact subset of $\C^* \setminus 2\D$ for every $\la \in X$;
\item[\rm (3)] $ S$ is upper semi-continuous;
\item[\rm (4)] if $C(\la)$ is the
component of  $\C^* \setminus { S}(\la)$ in $\C^*$ containing $2\D$,
for $\la \in X$, then $  C \cup { S}$ is upper semi-continuous and
$  C(\la)$ tends to $2\D$
as $|\la| \to 1$ in $X$.
\end{enumerate}

If $ S$ is a snare  we say that a function $s : \Delta \to \C^*$ is
{\em trapped by} $ S$ if $s(\la) \notin { S}(\la)$ for all $\la \in
X$.
\end{definition}

In (4), to say that $  C(\la)$ tends to $2\D$
as $|\la| \to 1$ in $X$ means the following.  For every $\eps > 0$ there exists $\de > 0$ such that $  C (\la) \subset (2+\eps)\D $ for all $\la\in X$ such that $|\la| > 1-\de$.  Note that $  C (\la)$ always contains $2\D$.

\newtheorem{Snare}[theorem]{Snare Lemma}
\begin{Snare} \label{snare}
Let ${ S}$ be a snare with domain $X$,
let $s: \Delta \to \C^*$ be a continuous function
trapped by $ S$ and suppose that,
for some $\la_0 \in X, \ |s(\la_0)| < 2.$ Then $|s(\la)| \leq 2$ for all
$\la \in \T$.
\end{Snare}
\proof
For $\la \in X$
let  $  C(\la)$ be the
component of  $\C^* \setminus { S}(\la)$ containing $2\D$.
Note that $  C(\la) \cup S(\la)$ is closed in $\C^*$, for 
if $w$ lies
outside this set then the component of $w$ in $\C^* \setminus S(\la)$
is a neighbourhood of $w$ disjoint from $  C(\la) \cup S(\la)$.
Let
$$
E = \{\la \in X : s(\la) \in \mathcal{C}(\la) \}.
$$
We shall prove that $E$ is open and closed in $X$.  Consider $z \in E$,
so that $s(z) \in {  C}(z)$.  Pick a connected open neighbourhood $V$
of $s(z)$ that meets $2\D$ and is such that
$\overline V \subset {  C}(z)$.
Then $\C^* \setminus \overline V$ is an open neighbourhood of ${\mathcal 
S}(z)$, and
so by the upper semi-continuity of $S$ there is a neighbourhood
$U$ of $z$ in $X$ such that
${S}(\la) \cap \overline V = \emptyset$ for all $\la \in
U$.  Now $U \cap s^{-1}(V)$ is a neighbourhood of $z$.  For  $\la \in
U \cap s^{-1}(V)$ the sets ${S}(\la)$ and $V \cup 2\D$ are disjoint,
and the latter set is a connected open set containing $2\D$.  Thus
$$
s(\la) \in V \cup 2\D \subset {  C}(\la).
$$
Hence $\la \in E$ for all $\la \in U\cap  s^{-1}(V)$, and so $E$ is open.

To show that $E$ is closed we consider the set
$$
  X \setminus E = \{ \la: s(\la) \notin {  C}(\la) \}.
$$
Consider $z \in X \setminus E$, so that $s(z) \notin {  C}(z).$
By the hypothesis that $s$ is trapped by $S$,
$s(z) \notin {  C \cup S}(z).$
Pick a closed neighbourhood $V$ of $s(z)$
disjoint from ${  C}(z) \cup {S}(z)$.
By the upper semi-continuity of $  C \cup S$
there is a neighbourhood $U$ of $z$ in $X$ such that ${  C}(\la) \subset
\C^* \setminus  V$ for all $\la \in U$.
For $\la \in U \cap s^{-1}(V)$ we have $s(\la) \in V$ and ${  C}(\la)
\subset \C^* \setminus  V$, so that $s(\la) \notin {  C}(\la)$.
Thus $U \cap s^{-1}(V)$ is a neighbourhood of $z$ contained in $ X 
\setminus E$.
Hence $X \setminus E$ is open, and so $E$ is closed in $X$.

Since $\la_0 \in X$ and $|s(\la_0)| < 2$ we have $s(\la_0) \in
{  C}(\la_0)$, and hence $\la_0 \in E$.  Thus $E$ is a non-empty subset
of $X$.  Since $E$ is open and closed in $X$, we must have $E=X$,
that is, $s(\la) \in {  C}(\la)$ for all $\la \in X$.

We can now deduce that $|s(\la)| \le 2$ for all $\la \in \T$.  For suppose
that there exists $\la_1 \in \T$ such that
$$
|s(\la_1)| = 2 + \varepsilon > 2
$$
for some $\varepsilon > 0$.  By the continuity of $s$ there exists
$\delta_1 > 0$ such that $|s(\la)| > 2+\half \varepsilon$ whenever
$|\la - \la_1| < \delta_1$, and by the fact that $S$ is a snare,
there exists $\delta_2 > 0$ such that
${  C}(\la) \subset (2+\half \varepsilon) \D$ whenever $\la \in X$ and
$1-\delta_2 < |\la| < 1$.  Since the closure of $X$ contains $\T$
we may find $\la \in X$ such that
$$
  |\la - \la_1| < \min \{ \delta_1, \delta_2 \}.
$$
Then $|s(\la)| > 2+ \half \varepsilon
$ but ${  C}(\la) \subset (2+\half \varepsilon) \D$,
contradicting the fact that $s(\la) \in   C(\la)$.  Hence $|s|\leq 2$ on $\T$.
\hfill \qed

\begin{lemma}\label{locus}
Let $\lambda \in \mathbb{D}$ and let $\sigma_1, \sigma_2,
\eta_1, \eta_2 \in \mathbb{T}$ be such that $\sigma_1 \neq \sigma_2,
\eta_1 \neq \eta_2$. Let
$$
  B= \{ \up(\la): \up \in \mathrm{Aut}~ \D,\  \up(\sigma_1) = \eta_1,
\ \up(\sigma_2) = \eta_2 \}.
$$
Then $ B$ is the intersection with  $\D$
of either a circle through $\eta_1$ and  $\eta_2$
or a straight line through $\eta_1$ and  $\eta_2$.
\end{lemma}
\proof Let
\begin{equation}\label{defchi}
\chi =
\frac{\sigma_2}{\sigma_1}\ \ \frac{\sigma_1-\lambda}{\overline \sigma_1-
\overline \lambda}\ \ \frac{\overline \sigma_2-
\overline \lambda}{ \sigma_2- \lambda}.
\end{equation}
A routine calculation with cross-ratios establishes the following.

If $\chi\eta_1 = \eta_2$ then $  B$ is the intersection with  $\D$
of the straight line through $\eta_1$ and  $\eta_2$.

If $\chi\eta_1 \neq \eta_2$ then $  B$ is the intersection with  $\D$
of the circle with centre
$$
\eta_1  \eta_2 \frac{\chi - 1}{\chi\eta_1 - \eta_2} \mbox{  and radius  }
\left| \frac{\eta_1 - \eta_2}{\chi\eta_1 - \eta_2} \right| ;
$$
moreover, this circle passes through $\eta_1$ and  $\eta_2$.
\qed

\section{A bound for $s$} \label{boundary}

In this section the main result, Lemma \ref{3sigmas}, provides an essential boundedness property ($|s|\leq 2$) for the candidate $h=(s,p):\D\to\C^2$ constructed on the assumption that Problem $\dia$ has a solution $p$.

 Typically there are $3$ boundary interpolation conditions (\ref{def2p}) in Problem $\dia$.  If we perform $3$ Schur reduction steps we can transform Problem $\diamondsuit$ to the search for a M\"obius function that maps the $\tau_\ell$ to $3$ given points on $\T$.  Such a M\"obius function exists if and only if the $3$ target points are in the same cyclic order as the $\tau_\ell$; this suggests that Problem $\diamondsuit$ can only be solvable by virtue of special properties of the $\tau_\ell$ and $m$. 

For brevity we introduce terminology for the hypotheses of Problem $\dia$.\\
{\bf Assumption A$\dia$}  Data $\la_j,\ s_j,\ p_j, \  j=1,2,3,$ satisfy
condition $\mathcal {C}_1$ extremally, with auxiliary extremal $m$, the function $q$  is the unique  Blaschke product of degree at most $2$ such that $q(\la_j)=\Phi(m(\la_j), s_j,p_j), \ j=1,2, 3,$ and the points $\tau_\ell, \  \ell=1,\dots, d(mq)$, are the distinct roots of the equation $mq(\tau)=1$.

\begin{lemma}\label{vals}
Under Assumption A$\dia$, let $p$ be a unimodular rational function and let $s$
be the rational function defined by equation {\rm (\ref{defs})}. 
If $p(\la_j)=p_j, \  j =1,2, 3,$
then, for any  M\" obius function $\up$ such that the
degree of $\varphi_\up$
is at most $3$,
$$
s(\lambda) \neq \frac{2}{\up (\lambda)}
$$
for every $\la\in\D$ such that $s(\la)^2\neq 4p(\la)$.
\end{lemma}

\begin{proof}
By hypothesis, condition
$\mathcal{C}_1$ holds, which is to say that
\begin{eqnarray*}
\lambda_j &\mapsto&  \varphi_\up(\lambda_j) ,\qquad  j =1,2, 3,
\end{eqnarray*}
are solvable Nevanlinna-Pick data for every $\up \in \B l_1$.
For $\up$ such that $\varphi_\up$ has degree at most $3$ it follows from Lemma
\ref{psiinner} that $\varphi_\up$ is a Blaschke product, and hence
has no poles in $\mathbb{D}$.  Suppose that, for some $\lambda \in
\mathbb{D},\ s(\lambda) = 2/ \up (\lambda)$.  Then $2-\up s$ vanishes at
$\lambda$, that is, $\up(\lambda) s(\lambda)=2$, and since $\lambda$ is not a pole of $\varphi_\up$,
the numerator $2 \up p -s$ must also vanish at $\lambda$.
 Hence
$s(\la)^2 - 4p(\la) = s(\la)^2 - 2\up(\lambda) s(\lambda)p(\la)= s(\la)(s-2\up p)(\la) =0$.
\end{proof}

The preceding lemma will give information about the values of $s$ whenever we can find $\up\in \Aut\D$ such that $d(\ph_\up)\leq 3.$   Observe that if $p$ satisfies equations (\ref{def2p}) then, by Lemma \ref{propsp},   $d(\ph_\up)\le 1+d(p)$.  However, by choosing $\up$ so that enough cancellations occur between the numerator and denominator $2\up p-s$ and $2-\up s$ of $\ph_\up$, we can arrange that $d(\ph_\up)=3$.

\begin{lemma} \label{3sigmas}
Let Assumption A$\dia$ hold, let $d(m)=1$ and let $p$ be a unimodular rational function of degree at most $4$ such that
\begin{eqnarray*}
p(\la_j) &=& p_j, \quad j=1,2,3,\\
p(\tau_\ell) &= & \bar m(\tau_\ell)^2, \quad  \ell=1,\dots, J,
\end{eqnarray*}
for some $J$ such that $d(q)+d(p)-3 \le  J \le 1+d(q).$  Suppose that $s_1^2\neq 4p_1$.
Let $s$ be the rational function defined by equation {\rm (\ref{defs})}.
If the relations 
\beq\label{defsig}
 (m^2p)(\si) =1, \qquad (mq)(\si) \neq 1
\eeq
have three distinct solutions
in $\T$ then $|s| \leq 2$ on $\T$.
\end{lemma}
\begin{proof}
Consider any pair $\sigma_i,\ \sigma_j$ of distinct points taken from the
three distinct solutions $\sigma_1,\ \sigma_2, \ \sigma_3 \in \T$ of (\ref{defsig}).
Denote by $\Upsilon_{ij}$ the set of M\"obius functions
$\up$ satisfying $\up(\sigma_i)=m(\sigma_i),\ \up(\sigma_j)=m(\sigma_j)$.
Apply Lemma \ref{propsp} with $\kappa =2$ and $n_1=n_2 =1$,
to show that
$\varphi_\up$
has degree at most
$$
d(p)+d(q)-J\le  3
$$
  for any $\up \in \Upsilon_{ij}$.
For any $\la \in \D$ let
$$
  B_{ij}(\la) =
\mathrm{clos}~ \{\up(\la): \up \in \Upsilon_{ij}  \} \subset \Delta.
$$
Note that $m \in \Upsilon_{ij}$, so that
$m(\la) \in   B_{ij}(\la)$.
By Lemma \ref{locus},
$  B_{ij}(\la)$ is the closed circular arc or straight line segment
in $\Delta$ joining
$ m(\sigma_i)$ to $ m(\sigma_j)$ and passing through $m(\la)$.

Let $D_1, D_2,$ and $D_3$ be pairwise disjoint closed circular discs 
contained in
$\Delta$, not containing $\la_1$ and such that $D_j$ is
  tangent to $\T$ at $\sigma_j$. Let $F$ be the set $\{\la\in\D:s(\la)^2=4p(\la)\}$; $F$ is finite by virtue of the identity (\ref{royal}).  Let
$X = \D \setminus \bigcup_j D_j \setminus F$, so that $X$ is a connected
open subset
of $\D$ whose closure contains $\T$. Since $s_1^2 \neq 4p_1$
and $(s,p)(\la_1)=(s_1,p_1)$ we have $\la_1 \in X$.

For $\la \in X$ with $|\la|$ close
to $1$, $m(\la)$ is not in the disc $m(D_j)$ touching $\T$ at $m(\sigma_j)$
and so $  B_{ij}(\la)$, which is the circular arc joining
$ m(\sigma_i)$ to $ m(\sigma_j)$ through $m(\la)$, is close to one of
the two arcs of $\T$ joining $ m(\sigma_i)$ to $ m(\sigma_j)$.  It follows
that, for every $\eps > 0$ there exists $\de >0$ such that 
\beq\label{ctyBij}
B_{ij}(\la) \subset \{z: 1-\eps \leq |z| \leq 1\}
\quad \mbox{ for all } \la\in X \mbox{ such that }|\la| > 1-\de.
\eeq

We claim that each $B_{ij}$  is upper semi-continuous on $X$.
Consider $\la_0\in X$ and a neighbourhood $U$ of $B_{ij}(\la_0)$.  There is a neighbourhood $V$ of $m(\la_0)$ such that the circular arc or straight line segment through $m(\si_i), \, m(\si_j)$ and any point in $V$ lies in $U$.  Then $B_{ij}(\la) \subset U$ for all $\la\in m^{-1}(V)$.  Thus $B_{ij}$ is upper semi-continuous.

Define a multifunction $  S$ from $X$  to $\C^*\setminus 2\D$ by
$$
  S_{ij}(\la) = \{ 2/z: z \in   B_{ij}(\la) \}
$$
and 
\[
S(\la) = \bigcup_{1\leq i<j\leq 3} S_{ij}(\la).
\]
$S(\la)$ is the union of the three circular arcs or straight line segments in $\C^*\setminus 2\D$ passing through $2/m(\la), 2\bar m(\si_i)$ and $2\bar m(\si_j)$ for $1\leq i<j\leq 3$ and is a compact subset of $\C^*\setminus 2\D$.  (For one $\la$ it will happen that $m(\la)=0$, but it does not matter.)
As in Definition \ref{defsnare}, let 
$ C(\la)$ be the connected component of $\C^*\setminus S(\la)$  in $\C^*$ containing $2\D$.   Since each $B_{ij}$ is upper semi-continuous, so are $S$ and $S\cup C$ on $X$.  As $\la\in X$ approaches the unit circle, $2/m(\la)$ approaches $2\T$ avoiding the three discs $2/m(D_j)\subset \C^*\setminus 2\D$ which are tangent to $2\T$ at $2\bar m(\si_j), \, j=1,2,3$.  It follows that, for any $\eps > 0$, there exists $\de>0$ such that $S(\la) \subset (2+\eps)\D$ whenever $\la\in X$ satisfies $|\la|> 1-|\de|$.  The same assertion holds with $S(\la)$ replaced by $C(\la)$, and so $C(\la)$ tends to $2\D$ as $|\la| \to 1$ in $X$.  Thus $  S$ is a snare on $X$.

We claim that $s$
is trapped by $  S$.  Consider any $\la \in X$.
We must show that $2/s(\la) \notin   B_{ij}(\la)$ for any pair of
indices $i,j$.  By Lemma \ref{vals}, $2/s(\la) \neq \up(\la)$ for any
$\up \in \Upsilon_{ij}$.
Notice that $  B_{ij}$ is defined as a closure, and so to conclude
that $2/s(\la) \notin   B_{ij}(\la)$
we must show that $s(\la) \neq 2\bar m(\sigma_i)$. Suppose the
contrary for some $\la \in X$.  Choose a sequence $(\up_n)$ in
$\Upsilon_{ij}$ converging
pointwise to $m(\sigma_i)$.  Since $\varphi_{\up_n}$ is inner we have
$$
\left | \frac{2\up_n(\la)p(\la)-2\bar m(\sigma_i)}
{2 - \up_n(\la)2\bar m(\sigma_i)} \right | = |\varphi_{\up_n}(\la)| < 1.
$$
Since the denominator tends to zero as $n \to \infty$, so does the
numerator and hence
$$
s(\la) = 2\bar m(\sigma_i), \quad p(\la) = \bar m(\sigma_i)^2.
$$
Thus $s(\la)^2= 4p(\la)$, contrary to choice of $\la \notin F$.
Hence $s(\la) \notin   S(\la)$, that is, $s$ is trapped
by $  S$.  

By equation (\ref{sj}),
$$
|s(\la_1)| = |s_1| < 2.
$$
By the Snare Lemma  $|s(\la)| \le 2$
for all $\la \in \T$.
\end{proof}
\begin{lemma}\label{Mgeq0}
Under Assumption A$\dia$,
\begin{equation}\label{pinter}
M \df \left[\frac{1-p_i \bar p_j}{1 - \lambda_i \overline
\lambda_j}\right]^3_{i, j=1} \ge 0
\end{equation}
and the Nevanlinna-Pick data $\la_j\mapsto p_j, \  j=1,2, 3$, are solvable.
\end{lemma}
\begin{proof}
Recall that by Proposition \ref{uppick}, Condition $\mathcal{C}_1$ is equivalent to the matrix inequality (\ref{pick}).
For any $\omega \in \mathbb{T}$ we can pick a sequence of
M\" obius functions $\upsilon$ converging uniformly on compact
subsets of $\mathbb{D}$ to the constant function $\omega$.  Take
limits along this sequence in the inequality (\ref{pick}) to infer that
$$
\left[\frac{1-p_i \bar p_j - \frac{1}{2} \omega(s_i - p_i\bar s_j)
-
\frac{1}{2} \overline \omega (\bar s_j - \bar p_j s_i)}{1-\lambda_i
\overline \lambda_j}\right]^3_{i, j=1} \ge 0.
$$
Integrate this inequality in $\omega$ with
respect to Lebesgue measure on $\T$ to obtain the Pick condition
(\ref{pinter}). \hfill $\Box$
\end{proof}

\section{Proof of Theorem \ref{main}}\label{mainproof}

 Let $\la_1,\la_2,\la_3$ be distinct points in $\D$ and let $(s_j,p_j)\in \G,\  j=1,2, 3$.  We must show that (1) the interpolation data 
\beq\label{3data}
\D\to \Ga: \la_j\mapsto (s_j,p_j), \quad j=1,2,3,
\eeq
are solvable by an aligned $\G$-inner function if and only if (2) the data satisfy condition $\mathcal{C}_1$ extremally and actively and Problem $\dia$ is solvable.

(2)$\Rightarrow$(1)
Suppose  that condition $\mathcal{C}_1$ holds extremally and actively
  and the corresponding Problem $\dia$ has a solution $p$.
 We shall construct a function $h$ in $\widetilde{\E_{1 3}}$
 such that $d(p)\le 4$ and $h(\la_j)=(s_j,p_j)$.  By \cite[Lemma 8.4]{ALY12}, $h(\D) \subset \G$ and thus, by Proposition \ref{aligned_nodes}, the function $h$ is aligned. 

Since condition $\mathcal{C}_1(\la, z)$ is active for the data $\la \mapsto z$ there exists  $m\in\aut\D$ such that the Nevanlinna-Pick data $\la_j\mapsto \Phi(m(\la_j),s_j,p_j)$, $\  j=1,2, 3$ are extremally solvable.  Hence there is a unique function $q$ in the Schur class that satisfies these interpolation conditions, and moreover $q$ is a Blaschke product of degree at most $2$.
  Let $\tau_\ell, \  \ell =1,\dots, 1+d(q),$ be the roots of the equation $mq(\tau) =1$.   Since $mq$ is a Blaschke product (of degree at most $3$),  each $\tau_\ell \in\T$.  

By hypothesis (2), the corresponding Problem $\dia$ has a solution $p$, that is, there exists a Blaschke product  $p$ of degree at most $4$  satisfying
\begin{equation} \label{solvedia}
\begin{array}{cclcl}
p(\la_j)&=&p_j, &\qquad&   j=1,2, 3, \\
p(\tau_\ell) &=& \bar m(\tau_\ell)^2,&  \qquad &\ell=1,\dots, d(mq).
\end{array}
 \end{equation}

  {\em We may assume that $s_1^2\neq 4p_1$}.  For if $s_j^2=4p_j$ for all $j$ then the Nevanlinna-Pick problem $\la_j\mapsto -\half s_j$ has a solution $f$ that is a Blaschke product of degree at most $2$.  Then $h=(2f, f^2)$ satisfies the interpolation conditions $\la_j\mapsto (s_j,p_j)$, and moreover every $m\in\aut\D$ is an auxiliary extremal.  Hence $h\in \widetilde{\E_{1 3}}$ and  condition (1)  of Theorem \ref{main} holds.  

{\em We may assume that the Pick matrix
\[
M = \left[ \frac{1-\bar p_i p_j}{1 - \bar \la_i\la_j}\right]_{i,j=1}^3 
\]
is positive definite}.
For by Lemma \ref{Mgeq0}, $M\geq 0$.  If $M$ is singular then, 
 by Proposition \ref{flat}, the interpolation data \eqref{3data} have a solution $h=(s,p)\in \hol(\D,\Ga)$. 
The function $\frac{2mp-s}{2-ms} \in\Schur$ satisfies the interpolation conditions $\la_j\mapsto \Phi(m(\la_j),s_j,p_j)$, $\  j=1,2, 3$, and so by uniqueness equals $q$; thus $h\in\widetilde{\mathcal{E}_{13}}$, as required.

Let $s$ be defined by equation \eqref{defs}, as in Lemma \ref{construct}, let $h=(s,p)$ and let $\ph_\ups=\frac{2\ups p-s}{2-\ups s}$ for any Blaschke product $\ups$.

   {\em We may assume that $d(p)=4$.}  It follows from the positivity of $M$ that $d(p)$ is either $3$ or $4$.
 Indeed, if
$$
p = c \prod_{i=1}^{d(p)} B_{\al_i},
$$
where $|c| = 1$, $\al_i \in \D$ and $d(p) \le 4$,
then
\[
1-\bar p_i p_j = 1 -\bar p(\la_i)p(\la_j)= \sum_{k=1}^{d(p)} \bar u_{ki}\left( 1-\bar B_{\al_k}(\la_i) B_{\al_k}(\la_j)\right) u_{kj}
\]
where 
\[
  u_{kj} = \prod_{1\leq \nu < k} B_{\al_\nu}(\la_j).
\]
Hence
\beq\label{exprM}
M=\left[ \frac{1-\bar p_i p_j}{1 - \bar\la_i \la_j}\right]_{i,j=1}^3=
\sum_{k=1}^{d(p)} (1-|\al_k|^2) v_k v_k^*
\eeq
for suitable column vectors $v_1, \dots , v_{d(p)}$.
By supposition the left hand side is a positive definite matrix,
and so $d(p) \ge 3$.   If $d(p)=3$ then $m^2p$ is a Blaschke product of degree $5$, and so the equation $m^2p=1$ has $5$ distinct roots, all in $\T$,
and the equation $mq=1$ has $1+d(q)$
distinct roots.  Thus the relations 
\beq\label{defsigbis}
m^2p(\si) =1, \qquad mq(\si)\neq 1
\eeq
have $4 - d(q) \ge 2$
solutions, which lie in $\T$. Let $\sigma$ be one of them.  Let $F$ be the finite set $\{\la\in\D: s(\la)^2=4p(\la)\}$ and
consider any $\lambda \in \D \setminus F$.  For any
$\mu \in \mathbb{D}$ we may choose a M\" obius function $\up $ such
that $\up (\sigma) = m(\sigma)$ and $\up (\lambda) =\mu$.
Indeed, $\up $ is given explicitly in terms of cross-ratios by
\beq\label{cross}
(\up (z), m(\sigma), \mu, 1/\overline \mu)  = (z, \sigma,
\lambda, 1/\overline \lambda).
\eeq
By Lemma \ref{cancels}, $\varphi_\up$ has a cancellation at $\si$ and so has degree at most $3$.   By Lemma \ref{vals}, $s(\lambda) \neq 2/\mu$.  Hence, by Lemma \ref{3sigmas},  $|s(\lambda)| \le 2$.  By continuity, the relation holds for
all $\la \in \Delta$.  By Lemma \ref{construct}, $(s,p)$ is the required interpolating function.  
Once again $\Phi\circ (s,p) \in\Schur$ satisfies the interpolation conditions \eqref{3data} and therefore, by uniqueness, equals $q$.  Hence $h\in\widetilde{\mathcal{E}_{13}}$.

{\em The relations \eqref{defsigbis} have $5-d(q)\geq 3$ distinct solutions in $\T$.}  For since $d(m)=1, \ d(p)=4$, the equation $m^2p=1$ has $6$ solutions, all in $\T$, which include the $1+d(q)$ solutions $\tau_\ell$ of the equation $mq=1$  by virtue of the equations (\ref{solvedia}).  Since $d(q) \leq 2$ we have $6-(1+d(q))\geq 3$ solutions of (\ref{defsigbis}).

Since  $p$ satisfies equations \eqref{solvedia}, we may apply Lemma \ref{propsp} with $J=d(mq)$.  Then the polynomial $\Pi$ of equation \eqref{defPi} is $1$ and so, by equation \eqref{exprs},  $s$ has the same denominator as $p$.  Thus $s$ is analytic in $\D$.    In Lemma \ref{3sigmas}  we may again take $J=1+d(q)$), and since there are $3$ distinct solutions of the relations \eqref{defsigbis}, we infer that $|s|\leq 2$ on $\T$.  By Lemma \ref{construct}, $h=(s,p)$ is an interpolating function in $\hol(\D,\Ga)$, and as before $\Phi\circ (m,h)=q$, and hence   $h\in \widetilde{\mathcal{E}_{13}}$.  Thus (1)$\Rightarrow$(2). \\

(1)$\Rightarrow$(2)
Suppose that the interpolation data \eqref{3data} are solvable by an aligned $\G$-inner function $h$.
By Proposition \ref{aligned_nodes}, $h \in  \widetilde{\E_{1 3}}$.
It means that the  function $h =(s,p) $ is rational, $p$ is  a Blaschke product of  degree less or equal $4$, and there exists $m \in \B l_{1}$ such that $d(m)=1$ and
$$
\frac{2 m p -s}{2- m s} \in \B l_{2}.
$$
By \cite[Proposition 5.1]{ALY12}, if $h=(s,p) \in \widetilde{\E_{1 3}}$, for the data 
$$ 
\la_j \to h(\lambda_j)= (s_j, p_j), \qquad j=1, 2, 3.
$$
condition $\mathcal{C}_1(\la, h(\la))$ is active. By Lemma \ref{pattau}, $p$ is a solution of the  corresponding
  Problem $\dia$. Thus condition (2) of  Theorem \ref{main} is satisfied. \hfill $\Box$

\begin{remark}  \rm If $h$ satisfies condition (1) of Theorem \ref{main} then, by Corollary \ref{aligned_extremal},  for {\em any} choice $\mu_1,\mu_2,\mu_3$ of distinct points in $\D$, the $3$-point interpolation data
\[
\mu_j\in\D \mapsto h(\mu_j)\in\G, \quad j=1,2,3,
\]
are extremally solvable. 
\end{remark}

\section{Caddywhompus functions} \label{caddy}
Theorem \ref{main} characterizes solvability of certain $3$-point interpolation problems by {\em aligned} $\G$-inner functions.  In this section we give examples of $\G$-inner functions of degree $4$ that are not aligned, and discuss their properties.

\begin{definition}\label{defcaddy}
A  rational $\G$-inner function $h=(s,p)$ is
{\em caddywhompus} if   $h(\D)\subset\G$,  the degree of $h$ is equal to $4$,  $h$ has at least $3$ distinct royal nodes in $\T$ and for every choice of $3$ distinct royal nodes $\omega_1,\omega_2,\omega_3$ in $\T$, the points $\half \overline{s(\omega_1 )}, \half \overline{s(\omega_2 )}, \half \overline{s(\omega_3 )} \in \T$  are not in the same cyclic order as $\omega_1,\omega_2,\omega_3$.
\end{definition}
Here we understand that  if one triple consists of distinct points and the other does not then the two triples do not have the same cyclic order (so that $(1,\ii, -1)$ and $(1,-1,1)$ do not have the same cyclic order).  The reason that cyclic orders play a role here is the following simple fact.  If $\la_1,\la_2,\la_3$ are distinct points in $\T$ and $\mu_1,\mu_2,\mu_3$ are any points in $\T$ then there exists $m\in\aut\D$ such that $m(\la_j) = \mu_j$ for each $j$ if and only if  the $\mu_j$ have the same cyclic order as the $\la_j$.

It follows from the definition that a rational $\G$-inner function $h \in  \hol(\D,\G)$ is caddywhompus if and only if it has degree $4$, has at least $3$ royal nodes in $\T$ and is not aligned.  Hence, by Proposition \ref{aligned_nodes}, $h\notin    \widetilde{\mathcal{E}_{13}}$.

There do exist both aligned and caddywhompus  $\G$-inner functions.

\begin{example} \label{excaddy}
\rm
(1) Consider again the degree $4$ $\G$-inner function $h$ of Example \ref{Examples5_2}(2).
  The royal nodes of $h$ in $\T$ are the three cube roots $\omega_j$ of $-1$, and $\half \overline{s(\omega_j)} = -\omega_j$ for each $j$.  Hence $h$ is aligned.\\

\noindent (2)  Let $0<\al<1$ and let $h$ be the symmetrization of the two Blaschke products $\la^2$ and $B_\al B_{-\al}$, that is,
$h(\la)= (\la^2+B_\al B_{-\al}(\la), \la^2 B_\al B_{-\al}(\la) )$.
The royal nodes of $h$ are the points $\la$ for which $\la^2=B_\al B_{-\al}(\la) = B_{\al^2}(\la^2)$, which are the points $\la=1, \ii,-1,-\ii$.  We may tabulate the royal nodes $\omega_j$ and the target values $\half \overline{s(\omega_j)}$: \\
\begin{center} \begin{tabular}{cccccc}
$ j$ & \quad & $1$ & $2$ & $3$ & $4$ \\
Royal node $\omega_j$ & \quad & $1$ & $\ii$ & $ -1$ & $-\ii$ \\
$\half \overline{s(\omega_j)}$ \quad & & $1$ & $-1$ & $1$ & $-1$.
\end{tabular} \end{center}
~\\

It is clear that, for any choice of $3$ royal nodes $\omega_j$, there are only $2$ corresponding target values $\half\overline{s(\omega_j)}$, and hence the target values are not in the same cyclic order as the nodes.  The degree $4\; \G$-inner function $h$ is therefore caddywhompus.\\

\noindent (3)   Let $-1 < \al < 1$ and let  $h$ be the symmetrization of the Blaschke products $\la^3$ and $B_\al$, so that 
\beq\label{greatex}
h(\la) = (\la^3+B_\al(\la), \la^3B_\al(\la)).
\eeq
Here
\[
(s^2-4p)(\la)=\frac{(\la^2-1)^2 (\al \la^2-\la+\al)^2}{(1-\al\la)^2}
\]
and so the royal nodes of $h$ are the points $1,\, -1$ and
\beq\label{rnodes}
\frac{1\pm\sqrt{1-4\al^2}}{2\al}.
\eeq
Thus if $|\al| < \half$ then $h$ has  $4$ royal nodes in $\R$, to wit $1, -1$  and the two points \eqref{rnodes}, of which one is in $\D$ and one lies outside $\Delta.$  When $\al = \pm \half$ the only royal nodes of $h$ are $1$ and $-1$.  Thus, for $|\al|\leq \half$, $h$ is neither aligned nor caddywhompus.   When $\half <|\al| < 1$, though, the nodes \eqref{rnodes} lie in $\T$, and so $h$ has four royal nodes in $\T$.   For example when $\al= -1/\sqrt{3}$ one has the royal node $\omega=\e^{\ii 5\pi/6}$ and $\half \overline{s(\omega)}=-\ii$.  The images of the nodes under $\half\bar s$ are in the opposite cyclic order to the nodes themselves. 
It follows that $\half\bar s$ maps every triple of royal nodes to a triple of distinct points in $\T$ in the opposite cyclic order.  Thus $h$ is caddywhompus.\\

\rm
\noindent (4)  Let $h(\la)=(\la^2+B_\al(\la), \la^2 B_\al(\la))$ where $-1<\al<1$.     The function $h$ is a $\G$-inner function of degree $3$ having $1$ as a royal node in $\T$.    There are $3$ cases.  If $\tfrac 13 <\al<1$ then
$h$ has $3$ distinct royal nodes in $\T$, to wit $1,\omega,\bar\omega$ where
\[
\omega =\frac{1}{2\al}(1-\al + \ii\sqrt{(3\al-1)(1+\al)}).
\]
Since $h$ has degree $3$ and has $2$ royal nodes $h$ is aligned.

For $\al\leq\tfrac 13$ there is only one royal node of $h$ in $\T$ (to wit, the point $1$), and so $h$ is not aligned.
When $-1<\al<\tfrac 13$ there are two other royal nodes, of which one is in $\D$ and the other is in $\C\setminus \Delta$.  When $\al=\tfrac 13$,
\[
(s^2-4p)(\la)= \frac{(\la-1)^6}{(3-\la)^2}
\]
and all the royal nodes coalesce at $1$.  Here  $h\in\mathcal{E}_{03}$ with the  auxiliary extremal  $m =1$ of degree $0$, and $h$ is $3$-extremal. 
\end{example}

The next result shows that if $3$-point interpolation data are generated by 
 localization of a caddywhompus function at $3$ points in $\D$ then the data do not satisfy $\mathcal{C}_1$ extremally  and actively.   
\begin{proposition}\label{notC1ext} 
Let $h=(s,p)$ be a caddywhompus $\G$-inner function and let $\la_1, \la_2, \la_3$ be distinct points in $\D$.
\begin{enumerate}
\item  The $\Ga$-interpolation data $\la_j \mapsto h(\la_j), \, j=1,2,3,$ do not satisfy condition $\mathcal{C}_1$ extremally and actively;
\item  if $s$ is injective on the set of royal nodes of $h$ in $\T$    then the  $\Ga$-interpolation data $\la_j \mapsto h(\la_j), \, j=1,2,3,$ do not satisfy condition $\mathcal{C}_1$ extremally. 
\end{enumerate}
\end{proposition}

\begin{proof} 
(1)  By Proposition \ref{aligned_nodes},  $h\notin \widetilde{\mathcal{E}_{13}}$.  By Definition \ref{defEclass} and Proposition \ref   {equivCnu}, 
condition $\mathcal{C}_1(\la,h(\la))$ does not hold extremally and actively.\\

\noindent (2)   We must show that there is no $m\in \B l_1$ such that the Nevanlinna-Pick data $\la_j\mapsto \Phi(m(\la_j), h(\la_j)), \;  j=1,2,3,$ are extremally solvable (Definition \ref{defCexly}).  Suppose there does exist such an $m$: then these Nevanlinna-Pick data are {\em uniquely} solvable, and the unique solution $q\in \B l_2$.  By \cite[Proposition 5.1]{ALY12}
\beq\label{yetagain}
q=\Phi\circ (m,h) =\frac{2mp-s}{2-ms}.
\eeq
If $d(m)=1$ then it follows that $h\in \widetilde{\mathcal{E}_{13}}$ and so, by Proposition \ref{aligned_nodes}, $h$ is aligned, a contradiction.  Alternatively, suppose that $m$ is a constant function.  Then since $d(q)\leq 2$ there must be at least two cancellations in equation \eqref{yetagain}, and hence it must be the case that $m=\half\overline{s(\omega)}$ for two distinct royal nodes $\omega$ of $h$ in $\T$, contrary to the hypothesis of injectivity.
Consequently the data do not satisfy $\mathcal{C}_1$ extremally.
\end{proof}
\begin{remark}\label{caddy3extrem}  
An example of a caddywhompus function that is $3$-extremal.  \rm
 In Proposition \ref{notC1ext}(2) we cannot delete the hypothesis of injectivity.  Let $h$ be the caddywhompus function in Example \ref{excaddy}(2).  Here we may choose $m$ to be the constant function $1$.  It is clear from the table that there are two cancellations in equation \eqref{yetagain}, at the royal nodes $1$ and $-1$, and so $\frac{2mp-s}{2-ms}$ has degree $2$.  Therefore $h\in\mathcal{E}_{03}$ and so is $3$-extremal.  Hence any $3$-point localization satisfies condition $\mathcal{C}_1$ extremally, and is therefore extremally solvable.
\end{remark}

On the other hand, Proposition \ref{notC1ext} tells us that if $h$ is the caddywhompus $\G$-inner function of Example \ref{excaddy}(3), as in equation \eqref{greatex} and $\la_1,\la_2,\la_3$ are any $3$ distinct points in $\D$, then the $\Ga$-interpolation data 
\beq\label{localize}
\la_j \mapsto h(\la_j), \; j=1,2,3,
\eeq
 do not satisfy $\mathcal{C}_1$ extremally.  
In fact the interpolation data in this example do not even satisfy condition $\mathcal{C}$ extremally, so that at present we have no way of showing that they are extremally solvable.  If they {\em are} extremally solvable then they constitute a counterexample to the Conjecture, for then, for some $r\in(0,1)$ close to $1$, the interpolation data $r\la_j\mapsto h(\la_j)$ satisfy $\mathcal{C}_1$ but are not solvable.  We therefore propose the following question.

{\em  Is every $\G$-inner rational function of degree $4$ having $4$ royal nodes in $\T$ $3$-extremal? }

An affirmative answer would refute our `$\Ga$-interpolation conjecture' \eqref{mainQ}.  More generally, if there is {\em any} $3$-extremal caddywhompus function  for which $s$ is injective on the set of royal nodes in $\T$  then the conjecture \eqref{mainQ} is false.

\section{Target data on the boundary}\label{dataOnBoundary}
Hitherto we have studied instances of Problem $I\Gamma$ in which the target points $z_j \in \Gamma$ lie in the {\em open} symmetrised bidisc $\G$.   For completeness this section discusses the case that some $z_j$ belongs to the topological boundary $\partial\Gamma$ of $\Gamma$.   The analysis of this case exhibits some interesting geometry of $\Gamma$.

In fact any map $h\in\hol(\D,\Gamma)$ satisfies either $h(\D)\subset \G$ or $h(\D) \cap \G=\emptyset$  (for example, \cite[Lemma 8.4]{ALY12}).  Thus a problem $I\Gamma$  can be solvable only if the target points are either all in $\G$ or all in $\partial\Gamma$, and consequently Problem $I\Gamma$ naturally splits into the problems $I\G$ and $I(\partial\Gamma)$.

Since $\partial\Gamma$ contains the embedded analytic disc
\[
D_\omega \df\{(\omega\la+\bar\omega,\la): \la\in\Delta\}
\]
for any $\omega\in\T$, one can easily write down examples of Problem $I(\partial\Gamma)$ which have non-constant solutions.
\begin{example}\label{superfic}
Let $\omega\in\T$, let $p_1,\dots,p_n\in\Delta$ and let $\la_1,\dots,\la_n$ be distinct points in $\D$.
The interpolation data
\beq\label{omw}
\la_j\in \D \mapsto (\omega p_j+\bar\omega,p_j) \in\Gamma, \quad j=1,\dots,n,
\eeq
are solvable if and only if the Nevanlinna-Pick data $\la_j\mapsto p_j$ are solvable, and in this case the solutions of the problem \eqref{omw} are $(\omega f+\bar\omega,f)$, where $f\in\Schur$ satisfies $f(\la_j)=p_j, \ j=1,\dots,n$.
\end{example}

Target data points in the {\em distinguished} boundary $b\Gamma$ of $\Gamma$ are special.  Recall \cite[Theorem 2.4]{AY04} that
\[
b\Gamma=\{(\la+\mu,\la\mu):\la,\mu\in\T\},
\]
whereas
\[
\partial \Gamma=\{(\la+\mu,\la\mu):\la\in\T, \mu\in\Delta\}.
\]
Note that $|\Phi_\omega(z)|=1$ for all $\omega\in\T$ and $z\in b\Gamma$.
\begin{lemma}\label{bGh}
If $h\in\hol(\D,\Gamma)$ and $h(\D)$ meets $b\Gamma$ then $h$ is constant.
\end{lemma}
\begin{proof}
Let $h\in\hol(\D,\Gamma)$ and let $h(\mu)\in b\Gamma$ for some $\mu\in\D$.   Then,  for all $\omega\in\T$,  $\Phi_\omega\circ h\in\Schur$ and 
$|\Phi_\omega\circ h(\mu)|=1$.
By the maximum principle $\Phi_\omega\circ h$ is constant on $\D$ for every $\omega\in\T$.   It is simple to deduce that $h$ is constant on $\D$.
\end{proof}
An immediate consequence of this lemma is that there are no non-trivial solvable interpolation problems $I(\partial\Gamma)$ in which a target data point lies in $b\Gamma$.
\begin{proposition}\label{bgamma}
Let $\la_1,\dots,\la_n$ be distinct points in $\D$, let $z_1\in b\Gamma$ and let $z_2,\dots,z_n\in\Gamma$.  The interpolation data
\[
\la_j\in\D \mapsto z_j \in\Gamma, \quad j=1,\dots,n,
\]
are solvable if and only if $z_1=\dots=z_n$, in which case the unique solution of the interpolation problem is the constant map $h(\la)=z_1$.
\end{proposition}
The case of target data lying in $\partial\Gamma\setminus b\Gamma$ is captured in Example \ref{superfic}.  To prove this we need some facts about $\Phi_\omega$.
\begin{lemma}\label{somefacts}
Let $\omega\in\T$ and $(s,p)\in\partial\Gamma\setminus b\Gamma$.  The following statements are equivalent.
\begin{enumerate}
\item $|\Phi_\omega(s,p)|=1$;
\item $\omega(s-\bar sp)=1-|p|^2$;
\item $(s,p) \in D_\omega$.
\end{enumerate}
\end{lemma}
\begin{proof}
The equivalence of (1) and (2) is \cite[Theorem 2.5]{AY04}.

Suppose (2).  We can write $s=\la+\mu, \, p=\la\mu$ for some $\la\in\T$ and $\mu\in\D$.  Then
\[
1=\frac{\omega(s-\bar sp)}{1-|p|^2} = \frac{\omega(\la+\mu-(\bar\la+\bar\mu)\la\mu)}{1-|\mu|^2} =\omega\la.
\]
Hence $\la=\bar\omega$ and $p=\bar\omega\mu$, and so $(s,p)=(\bar\omega+\omega p,p) \in D_\omega$.  Thus (2) implies (3).

Suppose (3): $s=\omega p+\bar\omega$ and $p\in\Delta$.  Then 
\[
\omega(s-\bar sp)= \omega(\omega p+\bar\omega -(\bar\omega\bar p+\omega)p)=1-|p|^2,
\]
and so (3) implies (2).
\end{proof}

\begin{proposition}\label{more superfic}
Let $\la_1,\dots,\la_n$ be distinct points in $\D$ and let $z_1,\dots,z_n\in\partial\Gamma\setminus b\Gamma$.  The interpolation data
\beq\label{superficdata}
\la_j\in\D \mapsto z_j \in\Gamma, \quad j=1,\dots,n,
\eeq
are solvable if and only if there exists $\omega\in\T$ and $p_1,\dots,p_n\in\Delta$ such that $z_j=(\omega p_j+\bar\omega,p_j)$ for $j=1,\dots,n$ and the Nevanlinna-Pick data
\[
\la_j\in\D \mapsto p_j\in\Delta, \quad j=1,\dots,n,
\]
are solvable.  In this case the solutions of the interpolation problem \eqref{superficdata} are the functions $(\omega f+\bar\omega,f)$, where $f\in\Schur$ satisfies $f(\la_j)=p_j, \ j=1,\dots,n$.
\end{proposition}
\begin{proof}
Sufficiency is Example \ref{superfic}.  To prove necessity, suppose that $h\in\hol(\D,\Gamma)$ is a solution of the problem \eqref{superficdata}.  Let $z_j=(s_j,p_j)$ for $j=1,\dots,n$.  Since $z_1\in\partial\Gamma\setminus b\Gamma$, a simple calculation shows that
\[
|s_1-\bar s_1p_1|= 1-|p_1|^2 > 0.
\]
There therefore exists a unique $\omega\in\T$ such that $\omega(s_1-\bar s_1p_1)= 1-|p_1|^2$.  It follows from Lemma \ref{somefacts}  that
$|\Phi_\omega(z_1)|=1$.  Now $\Phi\circ h \in\Schur$ satisfies
\[
\Phi_\omega\circ h(\la_1)= \Phi_\omega(z_1) \in\T,
\]
and so, by the maximum principle,  $\Phi_\omega\circ h$ is constant on $\D$.  Hence we have
\[
\Phi_\omega(z_j)=\Phi_\omega(z_1)\in\T, \quad \mbox{ for } j=2,\dots,n.
\]
Again by Lemma \ref{somefacts}, $z_j\in D_\omega$ for each $j$, that is
\[
z_j =(\omega p_j+\bar\omega, p_j)
\]
for $j=1,\dots,n$ and some $p_j\in\Delta$.  Furthermore, if $h=(s,p)$ then $p\in\Schur$ and $p(\la_j)=p_j$, and the Nevanlinna-Pick data $\la_j\mapsto p_j$ are solvable.
\end{proof}
\begin{remark} \rm
Each point $(s,p)$ of $\partial\Gamma\setminus b\Gamma$ lies in a unique disc $D_\omega$; the corresponding $\omega$ is given by $\bar\omega=(s-\bar sp)/(1-|p|^2)$.  Hence the condition in Proposition \ref{more superfic} that there exist $\omega\in\T$ such that $z_j=(\omega p_j+\bar\omega,p_j)$ for each $j$ can be written
\[
\frac{s_1-\bar s_1p_1}{1-|p_1|^2}= \dots = \frac{s_n-\bar s_n p_n}{1-|p_n|^2}.
\]

  Each pair of discs $D_\omega, \ D_\tau$, with $\tau\neq\omega\in\T$, intersects in the single point $(\bar\omega+\bar \tau, \bar\omega \bar\tau)$, which lies in $b\Gamma\setminus \{(2\omega,\omega^2): \omega\in\T\}$.  The point $(2\omega,\omega^2)$, on the `edge' of the M\"obius band $b\Gamma$ (see \cite[Theorem 2.4]{AY04}), lies on the unique disc $D_{\bar\omega}$.
\end{remark}

\section{Weak solvability does not imply solvability}\label{doesnot}

In this short section we justify the statement in Section \ref{duality} that weak solvability (recall Definition \ref {wkslvbl}) does not imply solvability for Problem $I\Gamma$.
For the proof which follows denote  by $H^2$ the Hardy Hilbert space on $\mathbb{D}$ and by $K$ the Szeg\H{o} kernel:
$$
K_\lambda (z) = K(z, \lambda) = (1-\overline \lambda z)^{-1},
\quad \lambda, z\in \mathbb{D}.
$$

\begin{proposition}\label{weakly-solv-not-solv} Let $n \ge 3$.
For any distinct points  $\la_1, \dots,\la_n$ in $\D$ there exist points 
$z_1, \dots, z_n$ in $\G$ such that the interpolation data 
$$
\la_j \mapsto z_j  \in \G,  \qquad j=1,\dots, n,
$$ 
are weakly solvable but not solvable.
\end{proposition}

\begin{proof} By  \cite[Theorem  12.4]{ALY12} there exist $z_1,\dots, z_n$ in $\G$ such that the 
 $\Gamma$-interpolation data 
\beq\label{again}
\lambda_j \mapsto z_j, \quad j=1,\dots,n,
\eeq
  satisfy $\mathcal{C}_{n-3}$ (and {\em a fortiori} $\mathcal{C}_0$) but not $\mathcal{C}_{n-2}$. Since $\mathcal{C}_{n-2}$ is necessary for solvability,  the $\Gamma$-interpolation data \eqref{again}
are unsolvable.

Let $z_j= (s_j, p_j)$, $j=1,\dots, n$. To say that the data satisfy $\mathcal{C}_0$ means that
\begin{equation}
\label{X_Cn-3}
  \|\Phi_\omega(S,P)\| \leq 1\quad \mbox{ for all } \omega\in\T
\end{equation}
where $S,P$ are the operators on 
\begin{equation}\label{defM-2}
\mathcal{M} = \mathrm{span}\ \{K_{\lambda_1}, \dots, K_{\lambda_n}\}  \subset H^2
\end{equation}
given by
\beq\label{defX-2}
SK_{\la_j}=  \bar s_j K_{\la_j}, \quad PK_{\la_j}=  \bar p_j K_{\la_j}, \quad  j=1,\dots,  n.
\eeq
Since $S$ has spectral radius $\max_j |s_j| < 2$, it follows from \cite[Theorem 1.2]{AY1} that $\si(S,P)\subset \G$ and $(S,P)$ is a $\Gamma$-contraction, which is to say that, for any $g\in\hol(\G,\D)$,
\[
\|g(S,P)\| \leq 1.
\]
By Pick's Theorem, as reformulated by Sarason \cite{sar}, the Nevanlinna-Pick data
\[
\la_j \mapsto g(s_j,p_j), \quad j=1,\dots,n
\]
are solvable, that is, the data $\la_j \mapsto z_j$ are weakly solvable.
\end{proof}

\section{More about extremally solvable data}\label{extsolvdata}

The purpose of this section is to show the relationship between four natural complex-geometric notions of extremal solvability and the notion we introduced in Section \ref{extremal}: extremal satisfaction of condition $\mathcal C$ is stronger than any of the geometric notions. 

Consider again the general interpolation data
\beq\label{verygen2}
\la_j\in D \mapsto z_j\in E, \qquad j=1,\dots,n,
\eeq
where $D$ is a domain and $E$ is a connected subset of $\C^N$ for some $N$.  We include the definition of extremal solvability (Definition \ref{extlySolv}) for the purpose of comparison.

\begin{definition} \label{exCoex}
The interpolation data \eqref{verygen2} are {\em extremally solvable} if they are solvable but there do not exist an open neighbourhood $U$ of the closure of $D$ and a map $h\in\hol(U,E)$ such that 
\beq\label{interpGen}
h(\la_j)=z_j \quad\mbox{ for } \quad  j=1,\dots,n.
\eeq
The interpolation data \eqref{verygen2} are {\em co-extremally solvable} if they are solvable but there do not exist a compact subset $K$ of the interior of $E$ and a map $h\in \hol(D,K)$ such that the conditions \eqref{interpGen} hold.

The data \eqref{verygen2} are {\em robustly solvable} if there is a neighbourhood $V_j$ of $\la_j$ in $D$ for $j=1,\dots,n$ such that, for all $\la_j'\in V_j$, the data 
\[
\la_j'\in D \mapsto z_j \in E, \qquad j=1,\dots,n,
\]
are solvable; otherwise the data \eqref{verygen2} are {\em barely solvable}. 

The data \eqref{verygen2} are {\em co-robustly solvable} if there is a neighbourhood $U_j$ of $z_j$ in $\C^N$ for $j=1,\dots,n$ such that $U_j \subset E$ and, for all $z_j'\in U_j$, the data 
\[
\la_j\in D \mapsto z_j' \in E, \qquad j=1,\dots,n,
\]
are solvable; otherwise the data \eqref{verygen2} are {\em co-barely solvable}. 
\end{definition}

\begin{remark}\label{holInvar}  \rm
{\rm (1)}  For distinct points $\la_1,\dots,\la_n\in\D$ let $\mathrm{Solv}_E(\la_1, \dots,\la_n)$ denote the set of points $(z_1,\dots,z_n)\in E^n$ such that the interpolation data 
\beq\label{Edata}
\la_j\in\D \mapsto z_j \in E, \qquad j=1,\dots,n,
\eeq
 are solvable, and by $\mathrm{Unsolv}_E(\la_1,\dots,\la_n)$ the complement of $\mathrm{Solv}_E(\la_1, \dots,\la_n)$ in $\C^{Nn}$.  Thus
\[
\mathrm{Solv}_E(\la_1, \dots,\la_n)= \{ (h(\la_1),\dots,h(\la_n)) : h\in\hol(\D,E) \}.
\]
Then the data \eqref{Edata} are co-barely solvable if and only if
\[
(z_1,\dots,z_n) \in \partial\mathrm{Solv}_E(\la_1,\dots,\la_n).
\]
{\rm (2)}
Robust solvability, co-robust solvability and co-extremal solvability are all holomorphically invariant: if $E$ is open and $\al:D\to D', \, \beta:E\to E'$ are biholomorphic maps then interpolation data $\la_j \mapsto z_j$ are robustly, co-robustly or co-extremally solvable for $\hol(D,E)$ if and only if the data $\al(\la_j) \mapsto \beta(z_j)$ are robustly, co-robustly or co-extremally solvable respectively for $\hol(D',E')$. 
  The analogous statement for extremal solvability is not true, since the isomorphism $\al$ does not necessarily extend to be analytic in a neighbourhood of the closure of $D$.  It follows that extremal solvability is not equivalent to co-extremal, bare or co-bare solvability in general.
\end{remark}
 
There is a simple implication between two of these notions of extremal solvability.
\begin{proposition}\label{coexbare}
Let $D$ be a bounded starlike domain and let $E$ be a domain.  If the interpolation data \eqref{verygen2} are robustly solvable then they are not co-extremally solvable.   If they are co-extremally solvable then they are barely solvable.
\end{proposition}
\begin{proof}
The second assertion is simply a restatement of the first.  We may suppose without loss of generality that $D$ is starlike about $0$.   Suppose the data \eqref{verygen2} are robustly solvable: then there exists $\eps > 0$ such that, for $1-\eps < r < 1$, the interpolation data 
\[
r\la_j \in D \mapsto z_j\in E
\]
are solvable.  Fix such an $r$ and let $\ph\in\hol(D,E)$ satisfy $\ph(r\la_j)=z_j, \ j=1,\dots,n$.   Let $\ph_r(\la)=\ph(r\la)$ for $\la\in D$.  Then $\ph_r$ is analytic in a neighbourhood of $D^-$.  Thus $\ph_r\in\hol(D,E)$, \ $\ph_r(D^-)$ is a compact subset of $E$ and $\ph_r(\la_j)=z_j$ for $j=1,\dots,n.$  Hence the data \eqref{verygen2} are not co-extremally solvable.
\end{proof}
There is a dual result to Proposition \ref{coexbare}, proved in much the same way.
\begin{proposition}\label{coexcobare}
Let $D$ be a domain and let $E$ be a  bounded domain with the property that, for $r\in(0,1)$, the closure of $rE$ is contained in $E$.  If the interpolation data \eqref{verygen2} are co-robustly solvable then they are not co-extremally solvable.   If they are co-extremally solvable then they are co-barely solvable.
\end{proposition}
\begin{remark} \rm
The property that the closure of $rE$ is contained in $E$ is strictly stronger than being starlike about $0$, as is shown by the example $r=\half$,
\[
E=\tfrac 14\D \cup(\D\cap \{z: \im z < 0\}) \subset \C.
\]
\end{remark}
\begin{proof}
Since the interpolation data \eqref{verygen2} are co-robustly solvable there exists $\eps>0$ such that $1<r<1+\eps$ implies that the data $\la_j\mapsto rz_j$ are solvable.  Pick any such $r$ and let $\ph\in\hol(D,E)$ satisfy $\ph(\la_j)=rz_j$ for $j=1,\dots,n$.  Let $\ph_r=\ph/r$.  Then $\ph_r\in\hol(D,E)$ and $\ph_r(\la_j)=z_j$.  The range of $\ph_r$ is contained in $r^{-1} E$, and so by hypothesis its closure is a compact set contained in $E$.  Hence the interpolation data $\la_j\mapsto z_j$ are not co-extremally solvable.
\end{proof}

  The extremal $\mathcal C$ condition, on which this paper is based, is stronger than all four geometric conditions.
\begin{theorem}\label{stronger4}
 Let   $\la_1, \dots,\la_n$  be distinct points in $\D$. If the  interpolation data
$$
\la_j \mapsto z_j  \in \G,  \quad j=1,\dots, n,
$$ 
are solvable and satisfy condition $\mathcal{C}$ extremally then the data are
\begin{enumerate}
\item   extremally solvable,
\item  co-extremally solvable,
\item  barely solvable and
\item  co-barely solvable.
\end{enumerate}
\end{theorem}
Some notation: for $\rho > 0$ and any point $(s,p)$ we define $\rho\cdot (s,p)=(\rho s,\rho^2 p)$.  We write $\rho\cdot \G$ for $\{\rho\cdot z: z\in\G\}$.\\
\begin{proof}  \noindent (1)   This is Theorem \ref{stronger}.   By hypothesis there exists $m\in\Schur$ such that the Nevanlinna-Pick data \eqref{active1} are extremally solvable; by Remark \ref{Cn} we may assume that $m\in \B l_n$.  Hence there exists $q \in \B l_{n -1}$ such that 
\begin{equation}
\label{Phiq-2}
\Phi(m(\la_j), z_j)= q(\la_j), \quad j =1, \dots, n.
\end{equation} 
   Suppose that the data $\la \mapsto z$ are not extremally solvable: there exists
an $r_0 > 1$ and a function $f\in \hol(r_0\D, \Gamma)$ such that $f(\la_j)=z_j$ for $j=1,\dots,n$. Since $f(\la_1)=z_1 \in \G$, $f(r_0\D)$ is  not contained in $\partial\Gamma$, and so, by \cite[Lemma 8.4]{ALY12},  $f(r_0\D) \subset \G$.

Pick any $r_1$ in the interval $(1,r_0)$: then $f(r_1\Delta)$ is a compact subset of $\G$.  Now 
\[
f(r_1\Delta) \subset \bigcup_{0<\rho<1} \,  \rho \cdot\G =\G,
\]
and hence there exists $\rho \in (0,1)$ such that $f(r_1\Delta)\subset \rho\cdot\G\subset\rho\cdot\Gamma$. 
Observe that, for $\la\in\Delta$ and $(s,p) \in\Gamma$, we have
\[
\Phi(\la,\rho\cdot(s,p)) =\Phi(\la,\rho s,\rho^2 p) = \frac{2\la\rho^2p-\rho s}{2-\la\rho s} = \rho\Phi(\rho\la, s, p) \in \rho\Delta.
\]
Thus
\[
\Phi(\Delta \times \rho\cdot\Gamma) \subset \rho\Delta \subset \D.
\]
Furthermore, $\Phi$ is analytic on $(\rho^{-1}\D) \times \rho\cdot\Gamma$.  Hence, by continuity of $\Phi$ and compactness of $\rho\cdot\Gamma$, there is a neighbourhood $U$ of $\Delta$ such that
\[
\Phi(U\times \rho\cdot\Gamma) \subset \D.
\]

Pick $r_2$ in the interval $(1,r_1)$ such that $m(r_2\D) \subset U$.  Then, for any 
$\la\in r_2\D \subset  r_1\D$, we have $m(\la) \in U$ and $f(\la)\in \rho\cdot\Gamma$,
and hence
\[
|\Phi(m(\la),f(\la))| < 1.
\]
Thus $\Phi\circ(m,f)$ belongs to the Schur class, and
\[
\Phi\circ(m,f)(\la_j) = \Phi\circ(m,h)(\la_j) =q(\la_j) \quad \mbox{ for } j=1,\dots,n.
\]
Hence $\Phi\circ(m,f)$ is a solution of the solvable Nevanlinna-Pick problem
\[
\la_j \mapsto q(\la_j), \quad j =1,\dots, n,
\]
as is $q\in \B l_{n-1}$. Any $n$-point Nevanlinna-Pick problem that is solved by an element of $\B l_{n-1}$ is  extremally solvable and has a unique solution, and so $\Phi\circ(m,f) = q$.
This yields a contradiction, since $\Phi\circ(m,f)$ maps $r_2\D$ into $\D$, whereas $q$ maps $r_2\D \setminus \Delta$ to the complement of $\Delta$.  
Thus the data $\la \mapsto z$ are extremally solvable, which is to say that (1) holds.

\noindent (2)  By hypothesis there exists $m\in\Schur$ such that the Nevanlinna-Pick data \eqref{active1} are extremally solvable.  Hence there exists $q \in \B l_{n -1}$ such that 
\begin{equation}
\label{Phiq2}
\Phi(m(\la_j), z_j)= q(\la_j), \quad j =1, \dots, n.
\end{equation} 
 Suppose that the data $\la \mapsto z$ are not co-extremally solvable; then there exists a compact subset $K$ of $\G$ and a function $h\in\hol(\D,K)$ such that $h(\la_j)=z_j$ for each $j$.  Since $\Phi\circ(m,h)$ is a solution of the extremally solvable Nevanlinna-Pick problem \eqref{active1}, we have $\Phi\circ(m,h)=q$.  On the other hand, since $|\Phi|<1$ on the compact set $\Delta \times K$ and $(m,h)(\D) \subset \Delta \times K$ we have
\[
\|q\|_\infty = \| \Phi\circ(m,h)\|_\infty \leq \sup_{\Delta \times K} |\Phi| < 1,
\]
which contradicts the fact that $q$ is a Blaschke product.  Hence the data are co-extremally solvable.\\

\noindent (3) follows from (2) and Proposition \ref{coexbare}.\\

\noindent (4)  follows from (2), Proposition \ref{coexcobare} and Lemma \ref{strongstarlike} below.
\end{proof}
\begin{lemma}\label{strongstarlike}
If $0<r<1$ then the closure of $r\G$ is contained in $\G$.
\end{lemma}
\begin{proof}
One can verify that the identity
\begin{align}\label{oldIdent}
|2&-zrs|^2-|2zrp-rs|^2 \notag \\
	&=r^2\left\{|2-zs|^2-|2zp-s|^2\right\}+4(1-r)(1+r-r\re(zs))
\end{align}
is valid for all $z\in\T, s,p\in\C$ and $r>0$ (this identity was used in \cite[page 380]{AY04}).  For $(s,p)\in \G$ we have $|s|<2$ and so
\[
1+r-r\re (zs) \geq 1-r.
\]
Moreover, by the inequality \eqref{Phicriter}, the first term on the right hand side of equation \eqref{oldIdent} is non-negative, and hence
\[
|2-zrs|^2-|2zrp-rs|^2 \geq 4(1-r)^2.
\]
On dividing through by $|2-zrs|^2$ we obtain
\[
1-|\Phi(z,rs,rp)|^2 \geq \frac{4(1-r)^2}{|2-zrs|^2}.
\]
Since $|2-zrs| \leq 2(1+r)$, it follows that
\[
1-|\Phi(z,rs,rp)|^2 \geq \frac{(1-r)^2}{(1+r)^2}
\]
for all $z\in\T, \ (s,p)\in\G$ and $0<r<1$.  By continuity, for fixed $r\in(0,1)$, any $(s,p)$ in the closure of $r\G$ and all $z\in\T$,
\[
1-|\Phi(z,s,p)|^2 \geq \frac{(1-r)^2}{(1+r)^2} >0
\]
and consequently $(s,p)\in \G$.
\end{proof}

\section{Concluding reflections}\label{conclude}

In this section we discuss the relevance of the main theorem and its method of proof to a problem that originally arose in control engineering.

 Nevanlinna-Pick interpolation theory has proved useful in control engineering: see for example \cite{Fr,DuPa}.  In the notation of this paper, it is Problem $IE$ where $E$ is the closed unit disc.
Nevanlinna-Pick theory is described in countless papers and books, including \cite{Pick,W,sar,bgr,AgMcC}.  The classical results extend with appropriate modifications to a very narrow class of other sets $E$, in particular to the case that $E$ is the closed unit ball of the space of $m\times n$ matrices.  For applications in engineering the theory would be much more useful if we could solve Problem $IE$ for a range of further sets $E$.  The simplest relevant non-classical target set appears to be $\Gamma$, and for this and other reasons many authors have studied the function theory of $\Gamma$.
 A summary, some background and references can be found in \cite{NJY11}.  It transpires that the theory is considerably more subtle than in the familiar classical cases, but is nevertheless amenable to analysis.

In this paper we study the $3$-point interpolation problem $I\Gamma$.  The attempt to reduce Problem $I\Gamma$ to a collection of classical Nevanlinna-Pick problems gave rise to the form of duality for $\G$ described in Section \ref{duality}.

As mentioned in Section \ref{duality} (equation \eqref{mainQ}) we have earlier conjectured that
{\em the $\Gamma$-interpolation data  
\[
\la_j\in\D \mapsto z_j \in\G, \quad j=1,\dots, n,
\]
  are  solvable  if and only if condition  $\mathcal{C}_{n-2}(\la, z)$ holds}.
The conjecture is true in the case $n=2$  \cite{AY04} but we still do not know if it holds when $n=3$.
Nevertheless our main result, Theorem \ref{main}, gives a strong partial result that can in principle be used to solve $3$-point interpolation problems $I\Gamma$ numerically, at least in a generic case.
The proof of the theorem describes a method of constructing aligned $\Gamma$-inner functions and thereby giving an approach to solving Problem $I\Gamma$ via a one-variable Nevanlinna-Pick interpolation problem, Problem $\dia$, for which a Pick-type solvability criterion is available (Corollary \ref{mainbis}).
However, the method will never yield a caddywhompus function, and so it is probably not fully general.

  Here is a high-level algorithm based on the proof of Theorem \ref{main}. 
Suppose given $3$-point interpolation data $\la_j\in\D \mapsto z_j\in\G$.  First test the necessary condition $\mathcal C_1$ for solvability given in Proposition \ref{uppick}; this entails checking the positivity of the pencil \eqref{pick} of $3\times 3$ matrices indexed by $\up\in \B l_1$.   Since $\B l_1$ is a compact set of $3$ real dimensions this should be numerically feasible. Consider first the case that condition $\mathcal C_1$ holds extremally.  Then, by Definition \ref{defCexly}, there is an auxiliary extremal $m\in\B l_1$ and a $q\in\B l_2$ with the properties described in Lemma \ref{construct}.  These Blaschke products can be found by a search over a low-dimensional compact set.  We anticipate that typically $m$ will have degree $1$, though there are cases in which $m$ is a constant.  Once $m$ and $q$ are known we may formulate the corresponding Problem $\dia$ (page \pageref{def1p}), which is a classical Nevanlinna-Pick problem, though with mixed interior and boundary interpolation conditions.  Problems of this type have been studied by numerous authors \cite{bolot, geo} and solvability criteria are as described in Corollary \ref{mainbis}.  If Problem $\dia$ is unsolvable then the initial Problem $I\Gamma$ is unsolvable.  If Problem $\dia$ has a solution $p$ then we may proceed as in the proof of Theorem \ref{main}.  Define
\[
s=2\frac{mp-q}{1-mq};
\]
then the Snare Lemma is used to prove that $|s|\le 2$ on $\T$ and hence that $(s,p)(\D)\subset\G$, and $(s,p)$ is the desired interpolating function in $\hol(\D,\G)$.

In the case that the interpolation data satisfy condition $\mathcal C_1$, but not extremally, one can choose $r\in(0,1)$ such that the data $r\la_j\mapsto z_j$ satisfy condition $\mathcal C_1$ extremally.  One may then proceed as above.  If the corresponding Problem $\dia$ is solvable then one can construct a solution $g$ of the modified interpolation problem; then the function $\la\mapsto g(r\la)$ is a solution of the initial problem.  However, if Problem $\dia$ is unsolvable, then we cannot conclude that the initial problem is unsolvable.  It may yet prove to be the case that Problem $\dia$ is {\em always} solvable -- if so, then the procedure we have outlined will in principle work provided that there is an auxiliary extremal $m$ of degree $1$.  In the exceptional case that condition $\mathcal C_1$ is inactive (that is, there are only {\em constant} auxiliary extremals $m$) we do not currently have a prescription.

The present results are only a first step towards a theory of interpolation that would meet the needs of control engineers.  Naturally one would like to solve interpolation problems with any number of nodes, and
it is natural to ask whether results about $\Gamma$ extend to the higher-dimensional symmetrised polydisc $\Gamma_N$.    D. Ogle found in his thesis \cite[Corollary 5.2.2]{ogle}  an analogue of the necessary condition $\mathcal{C}_0$ for interpolation into  $\Gamma_N$. However, when $N \ge 3$, this condition is insufficient for solvability even of two-point interpolation problems \cite[Observation 1.3]{bharali}.

JIM AGLER, Department of Mathematics, University of California at San Diego, CA \textup{92103}, USA\\

ZINAIDA A. LYKOVA,
School of Mathematics and Statistics, Newcastle University, Newcastle upon Tyne
 NE\textup{1} \textup{7}RU, U.K.~~\\

N. J. YOUNG, School of Mathematics and Statistics, Newcastle University, Newcastle upon Tyne NE1 7RU, U.K.
 {\em and} School of Mathematics, Leeds University,  Leeds LS2 9JT, U.K.

\begin{thebibliography}{99}\label{bibliog}

\bibitem{AgMcC} J. Agler and J. McCarthy, {\em Pick Interpolation and Hilbert Function Spaces}, Graduate studies in mathematics {\bf 44},  Amer. Math. Soc., Providence, R.I. 2002. 

\bibitem {ALY12}  J. Agler, Z. A. Lykova and N. J. Young, Extremal holomorphic maps and the symmetrized bidisc,  {\em  Proc. London Math. Soc.}  {\bf 106}(4) (2013)  781-818.

\bibitem{AY1} J. Agler and N. J. Young, A commutant lifting theorem for  a domain in ${\mathbb C}^2$ and spectral interpolation, {\em J.  Functional Analysis} {\bf 161} (1999) 452--477.

\bibitem {AY2} J. Agler and N. J. Young, Operators having the symmetrized bidisc as a spectral set, {\em Proc. Edinburgh Math. Soc.} {\bf 43} (2000) 195--210.


\bibitem{AY04T} J. Agler and N. J. Young,  The two-by-two
spectral Nevanlinna-Pick problem, {\em Trans. Amer. Math. Soc.}
 {\bf 356} (2004)  573--585.



\bibitem{AY04}  J. Agler and N. J. Young, The hyperbolic geometry of 
the symmetrized bidisc, {\em J. Geom. Anal.} {\bf 14} (2004) 375--403.


\bibitem{AY06}  J. Agler and N. J. Young, The complex geodesics of the symmetrized bidisc, {\em Inter. J. of Mathematics} {\bf 17}, no.4, (2006) 375--391.

\bibitem{magic}  J. Agler and N. J. Young,    The magic functions and automorphisms of a domain, {\em Complex Analysis and Operator Theory} {\bf  2} (2008) 383-404. 

\bibitem{bgr}  J. A. Ball, I. Gohberg and L. Rodman. \emph{Interpolation of Rational Matrix Functions}. Operator Theory: Advances and Applications \textbf{45} (Birkh\"auser Verlag, Basel, 1990).

\bibitem{bharali} G. Bharali, Some new observations on interpolation in the spectral unit ball, {\em Integral Equations and Operator Theory} 
{\bf 59}, no. 3, (2007)   329--343.

\bibitem{bolot}  {  V. Bolotnikov},  The boundary analog of the   Carath\'eodory-Schur interpolation problem, \emph{J. Approx. Theory}  \textbf{163} (4) (2011) 568--589.

\bibitem{BD} {  V. Bolotnikov and H. Dym},  On boundary interpolation for matrix valued Schur functions, \emph{AMS Memoirs} \textbf{181} (2006) Number 856 1--107.

\bibitem{ChenHu} G-N. Chen and Y-J. Hu, Multiple Nevanlinna-Pick interpolation with both interior and boundary data and its connection with the power moment problem, {\em Linear Algebra and its Applications} {\bf 323} (2001) 167 -194.

\bibitem{costara} C. Costara, The symmetrized bidisc and Lempert's theorem, {\em Bull. London Math. Soc.} {\bf 36} (2004) 656--662. 

\bibitem{Dieu} J. Dieudonn\'e,  {\em Foundations of Modern Analysis}, Academic Press, New York, 1960.


\bibitem{DuPa} G. Dullerud and F. Paganini, {\em A course in robust control theory: a convex approach}, Texts in Applied Mathematics {\bf 36}, Springer, 2000.

\bibitem{EZ}
A. Edigarian and W. Zwonek, Geometry of the symmetrised polydisc, {\em Archiv  Math.}, {\bf 84}  (2005) 364-374.

\bibitem{Fr} B. A. Francis, {\em A course in $H_\infty$ control theory},  Lecture Notes in Control and Information Sciences {\bf 88}, Springer Verlag Berlin, Heidelberg 1987.


\bibitem{geo}  {  D. R. Georgijevi\'c}, Mixed L$\ddot{\rm o}$wner and Nevanlinna-Pick Interpolation,  \emph{Integral Equations and Operator Theory}  \textbf{53} (2005) 247--267. 


\bibitem {JP} M. Jarnicki and P. Pflug, Invariant distances and metrics in complex analysis revisited, {\em Dissertationes Math. (Rozprawy Mat.)} {\bf 430} (2005) 1--192.

\bibitem{Ko98} S. Kobayashi, {\em Hyperbolic complex spaces}, Springer, New York, 1998.

\bibitem{Le86} L. Lempert, Complex geometry in convex domains, {\em Proc. Intern. Cong. Math.}, Berkeley, CA (1986) 759--765.

\bibitem{NiPfZw}
N. Nikolov, P. Pflug and W. Zwonek, The Lempert function of the symmetrized polydisc in higher dimensions is not a distance,
 {\em Proc. Amer. Math. Soc.} {\bf 135} (2007) 2921--2928.

\bibitem{ogle} D. Ogle,  Operator and Function Theory of the Symmetrized Polydisc, Ph. D. Thesis, Newcastle University, 1999,   http://hdl.handle.net/10443/1264 .


\bibitem{PZ} P. Pflug and  W. Zwonek, Description of all complex geodesics in the symmetrised bidisc, {\em Bull. London Math. Soc.} {\bf 37} (2005) 575--584.

\bibitem{Pick}  G. Pick.  \"Uber die Beschr\"ankungen analytischer Funktionen, welche durch vorgegebene  Funcktionswerte bewirkt werden, {\em Math. Ann.} {\bf 77} (1916) 7--23.


\bibitem{sar} D. Sarason, Generalized interpolation in $H^\infty$, {\em Trans. Amer. Math. Soc.} {\bf 127} (1967) 179--203.

\bibitem{Upmeier}  H. Upmeier, {\em  Toeplitz Operators and Index Theory in Several Complex Variables},
 Birkhauser Verlag, 1995.

\bibitem{W} J. L. Walsh, Interpolation and approximation
by rational  functions in the complex domain. Fourth edition. {\em American Mathematical Society Colloquium Publications},
Vol. XX, Amer. Math. Soc., Providence, R.I. 1965.

\bibitem{NJY11} N. J. Young,  Some analysable instances of $\mu$-synthesis. {\em Mathematical methods in systems, optimization and control}, Editors: H. Dym, M. de Oliveira, M. Putinar, Operator Theory: Advances and Applications, Vol. 222, Springer, Basel, 2012, pp. 349--366.
 
\end{thebibliography}
\end{document}